\documentclass[11pt]{amsart}
\usepackage{amssymb,pstricks,pstcol,pst-plot}

\definecolor{Cyan}{cmyk}{1,0,0,0}
\definecolor{GoldenRod}{cmyk}{0,0.10,0.84,0}
\definecolor{myblue}{rgb}{0.66,0.78,1.00}
\definecolor{Lavender}{cmyk}{0,0.48,0,0}
\definecolor{CornflowerBlue}{cmyk}{0.25,0.13,0,0}

\parskip=\smallskipamount

\newtheorem{theorem}{Theorem}[section]
\newtheorem{lemma}[theorem]{Lemma}
\newtheorem{corollary}[theorem]{Corollary}
\newtheorem{proposition}[theorem]{Proposition}

\theoremstyle{definition}

\newtheorem{example}[theorem]{Example}

\newtheorem{remark}[theorem]{Remark}

\newcommand{\C}{\mathbb{C}}

\newcommand{\N}{\mathbb{N}}

\newcommand{\Z}{\mathbb{Z}}
\renewcommand{\P}{\mathbb{P}}
\newcommand{\R}{\mathbb{R}}

\newcommand{\cC}{\mathcal{C}}

\newcommand{\cH}{\mathcal{H}}

\newcommand\wt{\widetilde}

\newcommand\hra{\hookrightarrow}

\newcommand\Hom{\mathrm{Hom}}

\newcommand\disc{\triangle}

\def\di{\partial}

\def\bs{\backslash}

\def\e{\epsilon}

\numberwithin{equation}{section}

\begin{document}
\title[Stein structures and holomorphic mappings]
{Stein structures and holomorphic mappings}
\author{Franc Forstneri\v c \& Marko Slapar}
\address{Institute of Mathematics, Physics and Mechanics, 
University of Ljubljana, Jadranska 19, 1000 Ljubljana, Slovenia}
\email{franc.forstneric@fmf.uni-lj.si, marko.slapar@fmf.uni-lj.si}
\thanks{Supported by grants P1-0291 and J1-6173, Republic of Slovenia.}
\thanks{The original publication is available at http://www.springerlink.com}
\thanks{http://dx.doi.org/10.1007/s00209-006-0093-0}

%
%
\date{July 5, 2006} 

\begin{abstract}
We prove that every continuous map from a Stein manifold $X$ to 
a complex manifold $Y$ can be made holomorphic by a homotopic 
deformation of both the map and the Stein structure on $X$. 
In the absence of topological obstructions 
the holomorphic map may be chosen to have pointwise maximal rank. 
The analogous result holds for any compact Hausdorff family of maps,
but it fails in general for a noncompact family.  
Our main results are actually proved for smooth almost complex
source manifolds $(X,J)$ with the correct handlebody structure. 
The paper contains another proof of Eliashberg's 
(Int J Math 1:29--46, 1990) homotopy characterization of 
Stein manifolds and a slightly different explanation of the construction
of exotic Stein surfaces due to Gompf 
(Ann Math 148 (2):619--693, 1998; J Symplectic Geom 3:565--587, 2005). 
\end{abstract}
\maketitle

\section{Introduction} 
A {\em Stein manifold} is a complex manifold which is 
biholomorphic to a closed complex submanifold of a Euclidean space 
$\C^N$ \cite{Ho}. The following is a simplified version of our main results, 
Theorems \ref{T6.1}, \ref{T6.3} and \ref{T7.1}.

%
%
%
%
\begin{theorem}
\label{Main1}
Let $X$ be a Stein manifold with the complex
structure operator $J$, and let $f\colon X\to Y$ be a continuous map 
to a complex manifold~$Y$. 
\begin{itemize}
\item[(i)] 
If $\dim_\C X \ne 2$, there exist a Stein complex structure $\wt J$ on $X$, homotopic to $J$,
and a $\wt J$-holomorphic map $\wt f\colon X\to Y$ homotopic to~$f$.
\item[(ii)] 
If $\dim_\C X =2$, there are an orientation preserving homeomorphism 
$h\colon X\to X'$ onto a Stein surface $X'$ and a holomorphic map 
$f'\colon X' \to Y$ such that the map $\wt f:= f'\circ h \colon X\to Y$ 
is homotopic to $f$.
\end{itemize}
The Stein structure $\wt J$ in (i), and the homeomorphism
$h$ in (ii), can be chosen the same for all members
of a compact Hausdorff family of maps $X\to Y$. 
\end{theorem}

More precisely, in case (i) we find a smooth homotopy $J_t\in{\rm End}_\R TX$ 
$(J_t^2=-Id,\ t\in[0,1])$ consisting of integrable (but not necessarily Stein)
complex structures on the underlying smooth manifold $X$, 
connecting the Stein structure $J_0=J$ 
with a new Stein structure $J_1=\wt J$, such that
$f$ is homotopic to a $\wt J$-holomorphic map. 
In case (ii) we get essentially the same statement after 
changing the smooth structure on $X$, i.e., 
the new Stein structure $\wt J$ on $X$ may be {\em exotic}.
More precise statements are given by Theorem \ref{T6.1}
for part (i), and by Theorem \ref{T7.1} for part (ii).

The question whether every continuous map from a Stein manifold
to a given complex manifold $Y$ is homotopic to a holomorphic map
is the central theme of the {\em Oka-Grauert theory}.
Classical results of Oka \cite{Oka}, Grauert \cite{G1,G2,G3}
and Gromov \cite{GOka} give an affirmative  answer 
when $Y$ is a complex homogeneous manifold or, 
more generally, if it admits a dominating spray. 
(See also \cite{FP1} and \cite{Fflex}.)
Recently this {\em Oka property} of $Y$ has been characterized in 
terms of a Runge approximation property for entire maps $\C^n\to Y$ 
on certain special compact convex subset of $\C^n$ \cite{FCAP,FFourier2}.
The Oka property holds only rarely as it implies in particular
that $Y$ is dominated by a complex Euclidean space, and this 
fails for any compact complex manifold of Kodaira general type.
For a discussion of this subject see \cite{Fflex}.
Although one cannot always find a holomorphic representative in each
homotopy class of maps $X\to Y$, Theorem \ref{Main1} gives 
a representative which is holomorphic with respect to 
some Stein structure on $X$ homotopic to the original one.

Even if the source complex manifold $X$ is not Stein, we can obtain a holomorphic 
map in a given homotopy class on a suitable Stein domain in $X$, 
provided that $X$ has a correct handlebody structure.

%
%
%
\begin{theorem}
\label{Main1bis}
Let $X$ be an $n$-dimensional complex manifold which admits a 
Morse exhaustion function $\rho\colon X\to\R$ without critical points of index $>n$. 
Let $f\colon X\to Y$ be  continuous map  to a complex manifold $Y$. 
\begin{itemize}
\item[(i)] If $n\ne 2$, there exist an open Stein domain $\Omega$ in $X$,
a diffeomorphism $h\colon X\to h(X)=\Omega$ which is diffeotopic to 
the identity map $id_X$ on $X$, and a holomorphic map $f'\colon \Omega\to Y$ such that 
$f'\circ h\colon X\to Y$ is homotopic to~$f$. 
\item[(ii)] 
If $n=2$, the conclusion in (i) still holds if $\rho$ has no critical points of 
index $>1$; in the presence of critical points of index $2$
the conclusion holds with $h$ a homeomorphism which is
homeotopic to $id_X$. 
\end{itemize}
\end{theorem}

Theorem \ref{Main1bis} immediately implies Theorem \ref{Main1}:
If $h_t\colon X\to h_t(X)\subset X$ $(t\in [0,1])$ is a diffeotopy
from $h_0=id_X$ to $h_1=h\colon X\to \Omega$ as in 
Theorem \ref{Main1bis} then $J_t:=h_t^*(J)$ is a homotopy
of integrable complex structures on $X$ connecting the original 
structure $J_0=J$ to a Stein structure $\wt J= h^*(J|_{T\Omega})$, 
and $\wt f=f'\circ h\colon X\to Y$ is a $\wt J$-holomorphic map
homotopic to $f$.

Our proof of Theorem \ref{Main1bis} shows that the only
essential obstruction in finding a holomorphic map in a given homotopy
class is that $X$ may be holomorphically `too large' to fit into $Y$.
This vague notion of `holomorphic rigidity' has several concrete manifestations,
for example, the distance decreasing property of holomorphic maps in most 
of the standard biholomorphically invariant metrics (such as Kobayashi's).
The problem can be avoided by restricting the size of the domain while
at the same time retaining the topological (and smooth in dimension $\ne 2$)
characteristics of~$X$.

The following simple example illustrates that Theorems \ref{Main1} 
and \ref{Main1bis} are optimal even for maps of Riemann surfaces.

\begin{example} 
\label{annuli}  
Let $X = A_r=\{z\in \C\colon 1/r <|z|<r\}$, 
and let $Y=A_R$ for another $R>1$. We have $[X,Y]=\Z$.
A homotopy class represented by an integer $k\in\Z$ admits a 
holomorphic representative if and only if $r^{|k|}\le R$,
and in this case a representative is $z \to z^k$. 
Since every complex structure on an annulus is biholomorphic 
to $A_r$ for some $r>1$, we see that at most finitely many homotopy 
classes of maps between any pair of annuli
contain a holomorphic map. 
The conclusion of Theorem \ref{Main1} can be obtained 
by a radial dilation, decreasing the value of $r>1$ 
to another value satisfying $r^k\le R$, which amounts 
to a homotopic change of the complex structure on $X$. 
This allows us to simultaneously deform any compact 
family of maps $X\to Y$ to a family of holomorphic maps, but it is 
impossible to do it for a sequence of maps belonging to 
infinitely many different homotopy classes. 
The problem disappears in the limit as $R=+\infty$ when 
$Y$ is the complex Lie group $\C^*=\C\bs \{0\}$ and
the Oka-Grauert principle applies \cite{G2}, \cite{Oka}.
The same phenomenon appears whenever the fundamental group 
$\pi_1(Y)$ contains an element $[\alpha]$ of infinite order 
such that the minimal Kobayashi length $l_N$ of loops in $Y$ 
representing the class $N [\alpha] \in\pi_1(Y)$ tends to 
$+\infty$ as $N\to +\infty$: A homotopically nontrivial 
loop $\gamma$ in $X$ with positive Kobayashi length $K_X(\gamma)$
can be mapped to the class $N[\alpha]$ by a holomorphic map 
$X\to Y$ only if $l_N \le K_X(\gamma)$, and this is possible 
for at most finitely many $N\in \N$. 
\end{example}

%
%
%
%

Our construction also gives holomorphic maps of maximal rank 
(immersions resp.\ submersions) provided that there are 
no topological obstructions. The following is a 
simplified version  of Theorem \ref{T6.3} below.

%
%
%
%
\begin{theorem}
\label{Main2}
Let $X$ be a Stein manifold of dimension $\dim X\ne 2$.
Assume that $f\colon X\to Y$ is a continuous map to a complex manifold $Y$
which is covered by a complex vector bundle map $\iota\colon TX\to f^*(TY)$ 
of fiberwise maximal rank. Then there is a Stein structure $\wt J$ 
on $X$, homotopic to $J$, and a $\wt J$-holomorphic map $\wt f\colon X\to Y$ of pointwise maximal 
rank which is homotopic to $f$. The analogous conclusion holds if $\dim X=2$ and 
$X$ admits a Morse exhaustion function $\rho\colon X\to \R$  without critical points of index $>1$.
\end{theorem}

Theorem \ref{Main2} is a holomorphic analogue 
of the Smale-Hirsch h-principle for smooth immersions \cite{Smale,Hirsch,Gbook} 
and of the Gromov-Phillips h-principle for smooth 
submersions \cite{Gfoliations,Ph}. 
The conclusion holds with a fixed Stein structure on $X$ provided that $Y$ 
satisfies a certain flexibility condition introduced (for submersions) 
in \cite{FFourier1}. For maps to Euclidean spaces see also 
\cite[\S 2.1.5]{Gbook} (for immersions) and \cite{FActa}  (for submersions).

%
%
%
%
An important source of Stein manifolds are the {\em holomorphically complete 
Riemann domains} $\pi \colon X \to\C^n$, $\pi$ being a 
locally biholomorphic map. These arise as the envelopes of 
holomorphy of domains in, or over, $\C^n$. Clearly every such manifold 
is holomorphically parallelizable, but the converse has been 
a long standing open problem:

\smallskip
{\em Does every $n$-dimensional Stein manifold $X$ with a trivial complex tangent 
bundle admit a locally biholomorphic map $\pi\colon X\to \C^n$~?}
\smallskip

In 1967 Gunning and Narasimhan gave a positive answer for 
open Riemann surfaces \cite{GN}. In 2003 the first author of this paper
proved that every parallelizable Stein manifold $X^n$ 
admits a holomorphic submersion $f\colon X\to \C^{n-1}$ \cite{FActa};
the remaining problem is to find a holomorphic function $g$
on $X$ whose restriction to each level set of $f$ has no critical points
(the map $(f,g)\colon X\to\C^n$ is then locally biholomorphic).
Theorem \ref{Main2} with $Y=\C^n$ $(n=\dim X)$ shows that
the above problem is solvable up to homotopy:

\begin{corollary}
If $(X,J)$ is a Stein manifold of dimension $n\ne 2$ whose 
holomorphic tangent bundle $TX$ is trivial then
there are a Stein structure $\wt J$ on $X$, homotopic to $J$, 
and a $\wt J$-holomorphic immersion $\pi\colon X\to\C^n$.
\end{corollary}

Note that every closed complex submanifold $X\subset \C^N$ 
with trivial complex normal bundle $T\C^N|_X/TX$ is parallelizable \cite{FR};
this holds in particular for any smooth complex  hypersurface in $\C^N$.

%
%
%
%

All our main results are actually proved in the class of smooth almost complex manifolds 
$(X,J)$ which admit a Morse exhaustion function $\rho\colon X\to \R$ without
critical points of index $> n= \frac{1}{2} \dim_\R X$
(see Theorems \ref{T6.1}, \ref{T6.3} and \ref{T7.1}). By Morse theory such $X$ is 
homotopically equivalent to a CW complex of dimension at most $n$ \cite{Mi}. 
Since the Morse indices of any strongly plurisubharmonic exhaustion function 
satisfy this index condition, this holds for every Stein manifold 
(Lefshetz \cite{Lef}, Andreotti and Fraenkel \cite{AF}, Milnor \cite{Mi}).
Conversely, if $(X,J)$ satisfies the above index condition 
and $\dim_\R X \ne 4$ then $J$ is homotopic to an integrable 
Stein structure on $X$ according to Eliashberg \cite{E}.
The present paper contains another proof 
of this important result, with an additional
argument provided in the {\em critical case}
when attaching handles of maximal real dimension 
$n=\frac{1}{2}\dim_\R X$. A detailed understanding of 
this construction is unavoidable for our purposes,
and assuming that the initial structure on $X$ 
is already integrable Stein (as in Theorem \ref{Main1}) 
does not really simplify our proof.

The story is even more interesting when $\dim_{\R} X=4$: 
{\em A smooth oriented four manifold without handles
of index $>2$ is homeomorphic to a Stein surface,   but the underlying  
smooth structure must be changed in general} (Gompf \cite{Go1}, \cite{Go2}).
Indeed, a closed orientable real surface $S$ smoothly embedded 
in a Stein surface $X$ (or in a compact K\"ahler surface 
with $b^+(X)>1$), with the only exception 
of a null-homologous $2$-sphere, satisfies the 
{\em generalized adjunction inequality}:
\begin{equation}
\label{adjunction}
		[S]^2 + |c_1(X)\cdotp S| \le -\chi(S).
\end{equation}
(See Chapter 11 in \cite{GompfStipsicz} and the papers 
\cite{FSteindomains,KM,LM,N,OSzabo}.)
For a 2-sphere the above inequality yields $[S]^2\le -2$. 
Taking $X=S^2\times \R^2=\C\P^1 \times\C$, the embedded 2-sphere 
$S^2\times\{0\} \subset X$ generates $H_2(X,\Z)=\Z$ and satisfies
$[S]^2=0$, hence $X$ does not admit any non-exotic Stein structure.

Nevertheless, {\em there is a bounded Stein domain in $\C^2$
homeomorphic to $S^2\times\R^2$}.  
This is a special case of  Gompf's result   
that for every tamely topologically embedded CW 2-complex 
$M$ in a complex surface $X$ there exists a topological 
isotopy of $X$ which is uniformly close to the identity on $X$
and which carries $M$ onto a complex $M'\subset X$ with 
a {\em Stein thickening}, i.e., an open Stein domain $\Omega\subset X$ 
homeomorphic to the interior of a handlebody with core $M$
\cite[Theorem 2.4]{Go2}.
In his proof, Gompf uses {\em kinky handles} of index 2 in 
each place where an embedded 2-handle with suitable properties cannot be found
in Eliashberg's construction. To obtain the correct manifold one must
perform an inductive procedure which cancels all superfluous loops 
caused by kinks, thereby creating {\em Casson handles} which are homeomorphic, 
but not diffeomorphic, to the standard index two handle $D^2\times D^2$ 
(Freedman \cite{Freedman}). In \S 7 we follow a similar path to construct a 
holomorphic map in the chosen homotopy class, performing the Casson tower 
construction simultaneously at a possibly increasing number of places.

\smallskip
{\em Organization of the paper.}
In \S 2 we recall the relevant notions from Stein and contact geometry. 
Sections \S 3 -- \S 5 contain preparatory lemmas. 
The main geometric ingredient is Lemma \ref{L3.1} which gives 
totally real discs attached 
from the exterior to a strongly pseudoconvex domain
along a complex tangential sphere. A main analytic ingredient 
is an approximation theorem for holomorphic maps to arbitrary complex manifolds 
(Theorem \ref{T4.1}). Lemma \ref{L5.1} provides an approximate extension 
of a holomorphic map to an attached handle. The main results are presented 
and proved in sections \S 6 
(for $\dim_\C X\ne 2$) and \S 7 (for $\dim_\C X =2$).

%
%
%
%
\section{Preliminaries}
%
%
%
%
We begin by recalling some basic notions of the handlebody theory;
see e.g.\ \cite{FQ}, \cite{GompfStipsicz}, \cite{Mi}. 
Let $X$ be a smooth compact $n$-manifold with boundary $\di X$,
and let $D^k$ denote the closed unit ball in $\R^k$.
A {\em $k$-handle} $H$ attached to $X$ is a copy of $D^k\times D^{n-k}$
smoothly attached to $\di X$ along $\di D^k\times D^{n-k}$,
with the corners smoothed, which gives a larger compact manifold with boundary. 
The central disc $D^k \times\{0\}^{n-k}$ is the {\em core} of $H$.
A {\em handle decomposition} of a smooth (open or closed) manifold 
$X$ is a representation of $X$ as an increasing union of compact domains
with boundary $X_j\subset X$ such that $X_{j+1}$ is obtained by
attaching a handle to $X_j$. (In the case of open manifolds one takes
the interior of the resulting handlebody.) By Morse theory
every smooth manifold admits a handlebody representation.

An {\em almost complex structure} on an even dimensional 
smooth manifold $X$ is a smooth endomorphism
$J\in \mathrm{End}_\R(TM)$ satisfying $J^2=-Id$.
The operator $J$ gives rise 
to the conjugate differential $d^c$, defined on functions by 
$\langle d^c\rho,v\rangle= - \langle d\rho, Jv\rangle$
for $v\in TX$, and the Levi form operator $dd^c$.
$J$ is said to be {\em integrable} if every point of $X$ 
admits an open neighborhood $U\subset X$ and a 
$J$-holomorphic coordinate map of maximal rank
$z=(z_1,\ldots,z_n) \colon U  \to \C^n$ $(n=\frac{1}{2}\dim_\R X$),
i.e., satisfying $dz\circ J=idz$; for a necessary and sufficient 
integrability condition see Newlander and Nirenberg \cite{NN}. 

If $h\colon X\to X'$ is a diffeomorphism and $J'$ is an
almost complex structure on $X'$, we denote by 
$J=h^*(J')$ the (unique) almost complex structure on $X$
satisfying $dh\circ J=J'\circ dh$; i.e., such that 
$h$ is a biholomorphism. Similarly we  denote 
by $J'=h_*(J)$ the push-forward of an almost complex
structure $J$ by $h$. A map $f'\colon X'\to Y$ to a 
complex manifold $Y$ is $J'$-holomorphic if and only
if $f=f'\circ h\colon X\to Y$ is $J$-holomorphic with 
$J=h^*(J')$.

An integrable structure $J$ on a smooth manifold $X$ 
is said to be {\em Stein} if $(X,J)$ is a Stein manifold;
this is the case if and only if there is a 
{\em strongly $J$-plurisubharmonic Morse exhaustion function}
$\rho\colon X\to \R$, i.e., $\langle dd^c \rho, v\wedge Jv \rangle >0$ 
for every $0\ne v\in TX$  (Grauert \cite{G4}).
The $(1,1)$-form $\omega = dd^c \rho=2i\di\overline \di \rho$ 
is then a symplectic form on $X$, defining a $J$-invariant Riemannian metric 
$g(v,w)= \langle \omega,  v\wedge Jw\rangle$ $(v,w\in TX)$. 
The Morse indices of such function $\rho$ are $\le n=\frac{1}{2} \dim_{\R} X$ 
and hence $X$ is the interior of a handlebody without
handles of index $>n$ \cite{AF,Mi}. 

A real subbundle $V$ of the tangent bundle $TX$
is said to be {\em $J$-real}, or {\em totally real}, if $V_x\cap JV_x=\{0\}$
for every $x\in X$; its complexification $V^\C=V\otimes_\R \C$ 
can be identified with the $J$-complex subbundle 
$V\oplus JV$ of $TX$. An immersion $G\colon D\to X$ 
of a smooth manifold $D$ into $X$ is $J$-real (or {\em totally real})
if  $dG_x(T_x D)$ is a $J$-real subspace of $T_{G(x)} X$ 
for every $x\in D$. 

Let $W$ be a relatively compact domain with 
smooth boundary $\Sigma=\di W$ in an almost complex manifold $(X,J)$. 
The set  $\xi = T\Sigma \cap J(T\Sigma)$ is a corank 
one $J$-complex linear subbundle of $T\Sigma$. 
Assume now that $\rho$ is a smooth function in a neighborhood of 
$\Sigma =\di W$ such that $\Sigma= \{\rho=0\}$, $d\rho \ne 0$ 
on $\Sigma$ and $\rho<0$ on $W$. Let $\eta:=d^c \rho|_{T\Sigma}$,
a one-form on $\Sigma$ with $\ker\eta= \xi$.
We say that $\Sigma$ is {\em strongly $J$-pseudoconvex},
or simply {\em $J$-convex}, if $\langle dd^c\rho, v\wedge Jv\rangle >0$
for all $0\ne v\in \xi$; this condition is independent 
of the choice of $\rho$. (We shall omit $J$ when it is clear which 
almost complex structure do we have in mind.) 
This implies that $\eta\wedge (d\eta)^{n-1}\ne 0$ on $\Sigma$ 
($n=\dim_\C X$) which means that $\eta$ is a {\em contact form} and 
$(\Sigma,\xi)$ is a {\em contact manifold} 
(see \cite[pp.\  338--340]{Du,E,Gbook}). 
A smooth function $\rho\colon X\to\R$ whose level
sets are $J$-convex outside of the critical points
is said to be $J$-convex.

An immersion $g\colon S\to \Sigma$ of a smooth manifold $S$
into a contact manifold $(\Sigma,\xi)$ is {\em Legendrian} 
if $dg(TS)\subset \xi$. In the case at hand,
when $\Sigma$ is the boundary of a strongly pseudoconvex domain,
another common expression is a {\em complex tangential immersion}.

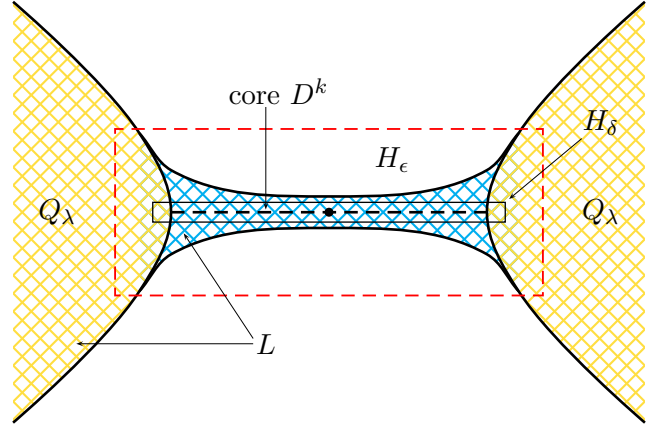
\begin{figure}[ht]
\psset{unit=0.7cm, xunit=1.2, linewidth=0.7pt} 
\begin{pspicture}(-5,-4)(5,4)

\pscustom[fillstyle=crosshatch,hatchcolor=Cyan, linestyle=none]     
{
\pscurve(-3,-1.5)(-2.5,-0.8)(0,-0.3)(2.5,-0.8)(3,-1.5) 
\psline[linestyle=dashed,linewidth=0.2pt](3,-1.5)(3,1.5)
\pscurve[liftpen=1](3,1.5)(2.5,0.8)(0,0.3)(-2.5,0.8)(-3,1.5)
\psline[linestyle=dashed,linewidth=0.2pt](-3,1.5)(-3,-1.5)
}

\pscustom[fillstyle=crosshatch,hatchcolor=GoldenRod,linewidth=1pt]
{\pscurve[liftpen=1](5,4)(3,1.5)(2.5,0)(3,-1.5)(5,-4)                 
}

\pscustom[fillstyle=crosshatch,hatchcolor=GoldenRod,linewidth=1pt]
{
\pscurve[liftpen=1](-5,4)(-3,1.5)(-2.5,0)(-3,-1.5)(-5,-4)             
}

\psline[linestyle=dashed,linewidth=1pt](-2.5,0)(2.5,0)                                 

\psecurve[linewidth=1pt](5,4)(3,1.5)(2.5,0.8)(0,0.3)(-2.5,0.8)(-3,1.5)(-5,4)           
\psecurve[linewidth=1pt](5,-4)(3,-1.5)(2.5,-0.8)(0,-0.3)(-2.5,-0.8)(-3,-1.5)(-5,-4)    

\psframe[linestyle=dashed, linecolor=red](3.4,1.6)(-3.4,-1.6)                          
\psframe[linewidth=0.5pt](2.8,0.2)(-2.8,-0.2)                                          

\psdot(0,0)

%
%
%
%
\rput(4.3,0){$Q_\lambda$}                                             
\rput(-4.3,0){$Q_\lambda$}
\rput(1,1){$H_\epsilon$}
\psline[linewidth=0.2pt]{->}(-1.2,-2.4)(-2.3,-0.4)
\psline[linewidth=0.2pt]{<-}(-4,-2.5)(-1.2,-2.5)
\rput(-1,-2.5){$L$}
\psline[linewidth=0.2pt]{->}(-1,2)(-1,0.05)
\rput(-0.8,2.3){core $D^k$}
\psline[linewidth=0.2pt]{<-}(2.85,0.25)(4,1.4)
\rput(4.3,1.7){$H_\delta$}

\end{pspicture}
\caption{A strongly pseudoconvex handlebody $L$}
\label{Fig1}
\end{figure}

Let $J_{st}$ denote the standard complex structure on $\C^n$. 
For a fixed $k\in \{1,\ldots,n\}$ let 
$z=(z_1,\ldots,z_n) = (x'+iy',x''+iy'')$, with $z_j=x_j+iy_j$,
denote the coordinates on $\C^n$
corresponding to the decomposition
\[
	\C^n=\C^k\oplus\C^{n-k}= 
	\R^k\oplus i\R^k \oplus\R^{n-k} \oplus i\R^{n-k}.
\]
Let $D=D^k\subset \R^k$ be the closed unit ball 
in $\R^k$ and $S=S^{k-1} =\di D$ its boundary sphere.
Identifying $D^k$ with its image in 
$\R^k \oplus \{0\}^{2n-k} \subset \C^n$ 
we obtain the core of the standard index $k$ handle 
\begin{equation}
\label{Hdelta}
	H_\delta = (1+\delta) D^k \times \delta D^{2n-k} \subset\C^n, \quad \delta>0.  
\end{equation}
A {\em standard handlebody of index $k$} in $\C^n$ 
is a set $K_{\lambda,\delta} = Q_\lambda \cup H_\delta$
for some $0<\lambda<1$ and $0<\delta<\frac{2\lambda}{1-\lambda}$ 
(Fig.\ \ref{Fig1}), where 
\begin{equation}  
\label{Qlambda}
	Q_\lambda = 
	\bigl\{ z=(x'+iy',z'') \in \C^{k}\oplus \C^{n-k}  \colon  
        |y'|^2  + |z''|^2 \le \lambda( |x'|^2-1) \bigr\}.                      
\end{equation}
The condition $\lambda<1$ insures that $Q_{\lambda}$ 
is strongly pseudoconvex, and the bound on $\delta$ implies
$(1+\delta)\di D^k\times \delta D^{2n-k}\subset Q_\lambda$. 

We shall need the following result of 
Eliashberg \cite[\S 3]{E}. (See also \cite{FK}.)

\begin{lemma}
\label{Eliashberg}
{\rm (Eliashberg)}
For every $\e>0$ and $\lambda>1$ there exist
a number $\delta\in (0,\epsilon)$ and a smoothly bounded,
strongly pseudoconvex handlebody $L\subset \C^n$ with core 
$Q_\lambda \cup D^k$ such that 
$K_{\lambda,\delta} \subset L \subset K_{\lambda,\epsilon}$ (Fig.\ \ref{Fig1}).
\end{lemma}

%
%
%
%

Eliashberg's construction gives an $L$ of the form
\[
	L= \bigl\{(x'+iy',z'') \in\C^n \colon |y'|^2 + |z''|^2 \le h(|x'|^2)  \bigr\}
\]	
where $h\colon [0,\infty]\to [\delta^2,\infty]$ is a smooth, increasing,
convex function chosen so that $L$ is a tube of constant radius 
$\delta$ around $D^k\subset \C^n$ over a slightly smaller ball 
$rD^k$ $(r<1)$, and $L$ equals $Q_\lambda$ over $r'D^k$
for some $r'>1$ close to~$1$.

We introduce the following (trivial) bundles over the disc 
$D \subset \R^k\oplus\{0\}^{2n-k}$:
\begin{eqnarray*}
    \nu\,' &=& \mathrm{Span} \Bigl\{\frac{\di}{\di y_1},\ldots, \frac{\di}{\di y_k}\Bigr\}\Big|_D = 
    D\times \bigl( \{0\}^k\oplus \R^k \oplus\{0\}^{2n-2k} \bigr), \\ 
     \nu\,''  &=&  \mathrm{Span} \Bigl\{\frac{\di}{\di x_j}, 
     		\frac{\di}{\di y_j} \colon  j=k+1,\ldots, n \Bigr\}\Big|_D 
            = D\times \bigl( \{0\}^{2k}\oplus \R^{2n-2k} \bigr), \\
     \nu    &=&  \nu\,'\oplus\nu\,'' = D\times \bigl(\{0\}^k\oplus \R^{2n-k}\bigr).
\end{eqnarray*}
Thus $\nu\,'=J_{st}(TD)$, $T^\C D = TD \oplus\nu\,'$,
and $T\C^n|_D = TD\oplus\nu = T^\C D\oplus \nu\,''$. 

Let $v\to S$ denote the (trivial) real line bundle over $S$
spanned by the vector field $\sum_{j=1}^k x_j \frac{\di}{\di x_j}$. 
Over $S$ we then have further decompositions
\[
	TD|_S = v \oplus TS, \quad 
	\nu\,'|_S =J_{st}(v) \oplus J_{st}(TS),  
        \quad TD|_S \oplus \nu\,'|_S \simeq v^\C \oplus T^\C S.
\]
Note that $T^\C S$ is a trivial complex vector bundle. 

Given a smooth embedding (or immersion) $G\colon D\to X$ 
of the disc $D=D^k \subset \C^n$ to a smooth $2n$-dimensional 
manifold $X$, a {\em normal framing} over $G$ is a homomorphism 
$\beta\colon \nu\to TX|_{G(D)}$ such that 
$
	dG_x \oplus \beta_x\colon T_x D\oplus \nu_x =T_x \C^n 
	\to T_{G(x)} X
$
is a linear isomorphism for every $x\in D$.

%
%
%
%
\section{Totally real discs attached to strongly pseudoconvex domains
along Legendrian spheres}
Let $W$ be an open, relatively compact domain with smooth strongly 
pseudo\-convex boundary $\Sigma=\di W$ in an almost complex manifold 
$(X,J)$ of real dimension $2n$. 
Choose a smooth defining function $\rho$ for $W$ which is
strongly $J$-plurisubharmonic near $\Sigma=\{\rho=0\}$.  
Let $w\subset TX|_\Sigma$ be the orthogonal complement of $T\Sigma$ 
with respect to the metric associated to the symplectic form $dd^c\rho$
(see \S 2); thus $w$ is spanned by the gradient of $\rho$ with respect 
to this metric. Then $Jw\subset T\Sigma$ and we have orthogonal 
decompositions $TX|_\Sigma = w\oplus T\Sigma= w\oplus Jw \oplus \xi$,
where $\xi=T\Sigma\cap J(T\Sigma)$.

Let $D=D^k$, $S=S^{k-1} = \di D$ and $v$ be as in \S 2.
An {\em embedding of a pair} $G\colon (D,S)\to (X\bs W,\Sigma)$
is a smooth embedding $G\colon D\hra X\bs W$ such that $G(S) = G(D)\cap \Sigma$
and $G$ is transverse to $\Sigma$ along $G(S)$. 
Such $G$ is said to be {\em normal to $\Sigma$} if $dG_x (v_x)=w_{G(x)}$ 
for every $x\in S$, i.e., $G$ maps the direction orthogonal 
to $S\subset \R^k$ into the direction orthogonal to 
$\Sigma\subset X$. The analogous definition applies to immersions.

The following lemma is a key geometric ingredient in the proof 
of all main results in this paper. Its proof closely follows the construction of a 
{\em special handle attaching triple} (HAT) in \S 2  of Eliashberg's paper 
\cite{E}, but with an additional argument in the critical case
$k=n\ne 2$ (see Remark \ref{framing} below). 
We thank Y.\  Eliashberg for his help in the proof of the
critical case (private communication, June 2005).

%
%
\begin{lemma}
\label{L3.1}
Let $W$ be an open, relatively compact domain with smooth 
strongly pseudoconvex boundary $\Sigma=\di W$ in an 
almost complex manifold $(X,J)$. Let $1\le k \le n=\frac{1}{2}\dim_\R X$, 
$D=D^k$, $S=\di D$. Given a smooth embedding 
$G_0\colon (D,S)\to (X\bs W,\Sigma)$, there is a regular 
homotopy of immersions $G_t\colon (D,S)\to (X\bs W,\Sigma)$ $(t\in [0,1])$ 
which is $\cC^0$ close to $G_0$ such that the immersion 
$G_1\colon D\to X\bs W$ is $J$-real and normal to $\Sigma$, 
and $g_1:=G_1|_S \colon S\hra \Sigma$ is a Legendrian embedding.
If $k<n$, or if $k=n\ne 2$, there exists an isotopy of embeddings 
$G_t$ with these properties. If $J$ is integrable 
in a neighborhood of $\Sigma\cup G_0(D)$ and $\Sigma$ 
is real-analytic then $G_1$ can be chosen real analytic.
\end{lemma}

As was pointed out in \cite[Note 2.4.2.]{E}, the topological obstruction 
in the case $k=n=2$ is essential. 
For example, there does not exist an embedded totally real 2-disc in $\C^2\bs B$,
attached to the ball $B\subset \C^2$ along a Legendrian curve in $\di B$, 
since by \cite{E} the resulting configuration would admit an open Stein
neighborhood diffeomorphic to $S^2\times\R^2$ in contradiction 
to the generalized adjunction inequality (\ref{adjunction}).

\begin{proof}  
The scheme of proof is illustrated on Fig.\ \ref{Fig2}.
First we find a regular homotopy from the initial disc $G_0\colon D\hra X\bs W$ 
to an immersed disc $G_1\colon D\to X\bs W$ which is attached with a correct normal framing
to $\di W$ along an embedded Legendrian sphere. 
Next we deform $G_1$ by a regular homotopy which is fixed near 
the boundary to a totally real immersed disc $G_2$, using
the h-principle for totally real immersions. Finally we show that, 
unless $k=n=2$, the construction can be done by isotopies of embeddings.

%
%
%
%
\begin{figure}[ht]
\psset{unit=0.8cm} 
\begin{pspicture}(0,-4)(12,4)

\pscircle[fillstyle=solid,fillcolor=myblue](3.5,0){2.5}           
\rput(3.5,0){$W$}						  

\psarc[linewidth=1.5pt,linecolor=green](7,0){2}{-135}{135}                              
\psecurve[linecolor=red,linewidth=1.5pt](4,1)(5.1,1.9)(5.6,2.3)(7,2)(8,2.5)             
\psecurve[linecolor=red,linewidth=1.5pt](4,-1)(5.1,-1.9)(5.6,-2.3)(7,-2)(8,-2.5)        
\psarc[linecolor=blue,linewidth=1.5pt](7,0){2.71}{-123}{123}                            

\psline[linewidth=0.2pt]{<-}(6.2,2.25)(7,3.1)                           
\psline[linewidth=0.2pt]{<-}(7.2,2.05)(7.2,3.1)
\rput(7.2,3.4){$G_1$}

\psline[linewidth=0.2pt]{->}(8.4,0)(8.95,0)                             
\rput(8,0){$G_0$}

\psline[linewidth=0.2pt]{<-}(9.75,0)(10.3,0)                            
\rput(10.7,0){$G_2$}

\psdots(5.12,1.91)(5.6,1.4)(5.12,-1.91)(5.6,-1.4)                       

\rput(4.9,1.7){$g_1$}                                                   
\rput(5.35,1.2){$g_0$}
\rput(4.9,-1.7){$g_1$}
\rput(5.35,-1.2){$g_0$}

\end{pspicture}
\caption{Deformations of an attached disc} 
\label{Fig2}
\end{figure}
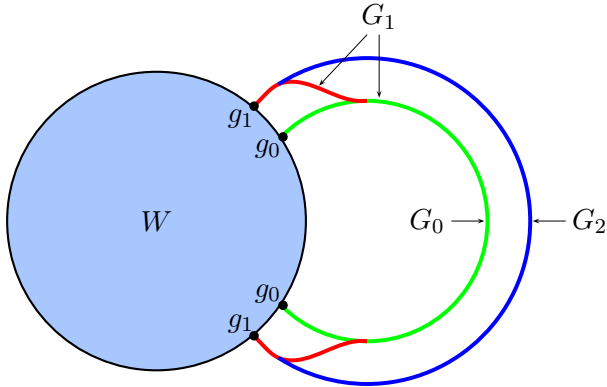

Set $g_0=G_0|_S\colon S\hra \di W$.
By a correction of $G_0$ along $S$ (keeping $g_0$ fixed) 
we may assume that it is normal to $\Sigma$,
i.e., such that $l_0:=dG_0|_v$ maps $v$ to $w|_{g_0(S)}$. 
Choose a complex vector bundle isomorphism
\[
	\phi_0\colon T\C^n|_D=D\times\C^n \to TX|_{G_0(D)},
	\quad \phi_0\circ J_{st} = J\circ \phi_0
\] 
covering $G_0$. We shall use the coordinates on $\C^n$ introduced
in \S 2. The vector field $\tau = \sum_{j=1}^k x_j \frac{\di}{\di x_j}$
is outer radial to the sphere $S=\di D$ in $\R^k\times\{0\}^{2n-k}$. 
Let $\wt \tau$ be the unique nonvanishing vector field on $\C^n$ over $S$ 
satisfying $\phi_0(\wt  \tau_x) = \ell_0(\tau_x)$ for every $x\in S$. 
By dimension reasons there exists a map 
$A\colon D\to GL_n(\C)$ satisfying $A_x \tau_x=\wt \tau_x$
for $x\in S$. Replacing $\phi_0$ by $\phi_0\circ A$ we may 
(and shall) assume from now on that $\phi_0|_v = \ell_0$. 
A further homotopic correction of $\phi_0$ insures that 
$\phi_0(T^\C S \oplus \nu''|_S)=\xi|_{g_0(S)}$,
thereby providing a trivialization of the latter bundle.

Write $\phi_0=\phi'_0\oplus \phi''_0$ where $\phi'_0=\phi_0|_{T^\C D}$
and $\phi''_0=\phi_0|_{\nu\,''}$ (we use the notation of \S 2). 
Setting $\psi_0:=\phi_0|_{T^\C S}$ we thus have 
\[
	\phi'_0|_{T^\C D|_S} = \ell_0^\C \oplus \psi_0 
	\colon v^\C \oplus T^\C S  \to TX|_{g_0(S)} = w^\C \oplus \xi|_{g_0(S)}.
\]
Note that $\psi_0\oplus \phi''_0\colon T^\C S \oplus \nu\,''|_S \to \xi|_{g_0(S)}$
is a complex vector bundle isomorphism.  
Furthermore, there is a homotopy of real vector bundle monomorphisms  
$\iota_s \colon TD\hra TX|_{G_0(D)}$ $(s\in [0,1])$
satisfying
\[
	\iota_0= dG_0, \quad \iota_1=\phi_0|_{TD}, \quad
	\iota_s|_v =\ell_0 \colon v\to w|_{g_0(S)} \ \ (s\in[0,1]).
\]

Consider the pair $(g_0,\psi_0)$ consisting of the embedding 
$g_0\colon S\hra \Sigma$ and the $\C$-linear embedding 
$\psi_0\colon T^\C S\hra \xi|_{g_0(S)}$
of the complexified tangent bundle of $S$ 
(a trivial complex vector bundle of rank $k-1$)
into the contact subbundle $\xi \subset T\Sigma$ over 
the map $g_0$.  By the {\em Legendrization theorem} 
of Gromov (\cite{Gbook}, p.\ 339, (B')) and Duchamp \cite{Du} 
there exists a Legendrian embedding $g_1\colon S\hra\Sigma$
whose complexified differential $\psi_1:=d^\C g_1$ 
is homotopic to $\psi_0$ by a homotopy of $\C$-linear 
vector bundle embeddings $\psi_t \colon T^\C S\hra \xi$ 
$(t\in [0,1])$.

Let $\Hom_{inj}(TS,T\Sigma)$ denote the space of all
fiberwise injective real vector bundle maps  
$TS\hra  T\Sigma$. Consider the path in $\Hom_{inj}(TS,T\Sigma)$
beginning at $dg_0$ and ending at $dg_1$, consisting of 
the homotopy $\iota_s|_{TS}$ $(s\in [0,1])$ followed by 
the homotopy $\psi_t|_{TS}$ $(t\in [0,1])$ (left and top side
of the square in Fig.\ \ref{Fig3}). By Hirsch's one parametric 
h-principle for immersions \cite{Gbook,Hirsch} this path 
can be deformed in the space $\Hom_{inj}(TS,T\Sigma)$ (with fixed ends) 
to a path $dg_t\colon TS\hra T\Sigma|_{g_t(S)}$ where 
$g_t\colon S\to \Sigma$ $(t\in [0,1])$ is a regular homotopy of 
immersions from $g_0$ to $g_1$.
We can insure that $\psi_t$ covers the base map $g_t$ 
for all $t\in [0,1]$. This gives a two parameter 
homotopy $\theta_{t,s}\in \Hom_{inj}(TS,T\Sigma)$ 
for $(t,s)\in [0,1]^2$ satisfying the following
conditions (Fig.\ \ref{Fig3}):
\begin{itemize}
\item[(i)]    $\theta_{t,0}= dg_t$ (bottom side),
\item[(ii)]   $\theta_{t,1}= \psi_t|_{TS}$ (top side),
\item[(iii)]  $\theta_{0,s}= \iota_s|_{TS}$ (left side; 
hence $\theta_{0,0}=dg_0$ and $\theta_{0,1}=\psi_0|_{TS}$), 
\item[(iv)]   $\theta_{1,s}= dg_1$ (right side), and 
\item[(v)]    $\theta_{t,s}$ covers $g_t$ for every $t,s\in [0,1]$.
\end{itemize}

%
%
%
%

\begin{figure}[ht]
\psset{unit=0.6cm, linewidth=0.7pt}  

\begin{pspicture}(-2,-1)(7,6)
\psframe[fillstyle=crosshatch,hatchcolor=GoldenRod](0,0)(5,5)

\psline[linewidth=2pt,linecolor=red](5,0)(5,5)

\psdots(0,0)(0,2.5)(0,5)(2.5,0)(2.5,2.5)(2.5,5)(5,0)(5,5)

\rput(0,-0.4){$dg_0$}
\rput(2.5,-0.4){$dg_t$}
\rput(5,-0.4){$dg_1$}

\rput(-1.2,2.5){$\theta_{0,s}=\iota_s$}
\rput(2.9,2.9){$\theta_{t,s}$} 
\psline[linewidth=0.2pt]{<-}(5.1,2.5)(5.8,2.5) 
\rput(7,2.5){constant}

\rput(0,5.5)  {$\psi_0$}
\rput(2.5,5.5){$\psi_t$}
\rput(5,  5.5){$\psi_1$}

\end{pspicture}
\caption{The homotopy $\theta_{t,s}$}
\label{Fig3}
\end{figure}
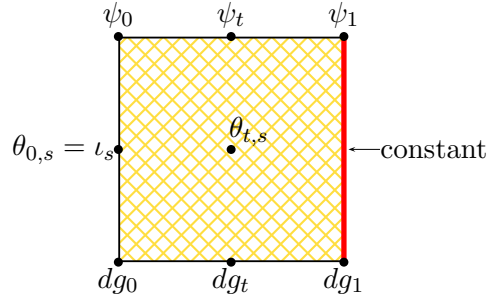

We can extend $g_t$ to a regular homotopy of immersions 
$G_t\colon (D,S)\to (X\bs W,\Sigma)$ 
$(t\in [0,1])$ which are normal to $\Sigma$, 
beginning at $t=0$ with the given map $G_0$. 
Let $\ell_t := dG_t|_v \colon v\to w|_{g_t(S)}$.
By the homotopy lifting theorem there exists a homotopy 
of $\C$-linear vector bundle isomorphisms $\phi_t$ covering $G_t$, 
\[
	\phi_t =\phi'_t\oplus \phi''_t \colon 
	T\C^n|_D = T^\C D \oplus \nu\,'' \to TX|_{G_t(D)},
	\quad t\in [0,1],
\] 
beginning at $t=0$ with the given map $\phi_0$, such that
over $S=\di D$ we have 
\[
	\phi'_t = \ell_t^\C \oplus \psi_t, \quad t \in [0,1],
\]
and  $dG_1=\phi_1$ on $TD|_S$.   

Set $\widetilde \theta_{t,s} = \ell_t\oplus \theta_{t,s}
\colon TD|_S \hra TX|_{g_t(S)}$ for $t,s\in [0,1]$
(a real vector bundle monomorphism over $g_t$).
From the above properties (i)--(v) of $\theta_{t,s}$ 
we obtain
\begin{itemize}
\item[(i')]    $\widetilde \theta_{t,0} = \ell_t\oplus dg_t= dG_t|_{TD|_S}$ (bottom side),
\item[(ii')]   $\widetilde \theta_{t,1} = \ell_t\oplus \psi_t|_{TS} =\phi_t|_{TD|_S}$ (top side),
\item[(iii')]  $\widetilde \theta_{0,s}= \iota_s|_{TD|_S}$ (left side), 
\item[(iv')]   $\widetilde \theta_{1,s}= \ell_1\oplus dg_1 = dG_1|_{TD|_S}$ (right side), and 
\item[(v')]    $\widetilde \theta_{t,s}$ covers $g_t$ for every $t,s\in [0,1]$.
\end{itemize}

We wish to extend the monomorphisms 
$\widetilde \theta_{t,s} \colon TD|_S \hra TX|_{g_t(S)}$ 
to real vector bundle monomorphisms 
$\Theta_{t,s}\colon TD\to TX$ ($t,s \in [0,1]$) covering the 
immersions $G_t \colon D\hra X$. Such extension already exists for $(t,s)$ 
in the bottom, top and left face of the parameter square $[0,1]^2$ 
where we respectively take $dG_t$, $\phi_t|_{TD}$ and $\iota_s$ 
(properties (i'), (ii') and (iii')). 
The homotopy lifting property provides an extension 
$\Theta_{t,s}$ for all $(t,s)\in [0,1]^2$,
with the given boundary values on the bottom, top and left side 
of $[0,1]^2$. (See Fig.\ \ref{Fig4}; the front and the back face  
belong to the homotopy $\wt\theta_{t,s}$ over $S=\di D$,
compare with Fig.\ \ref{Fig3}.)
Over the right face $\{t=1\}$ we thus obtain a homotopy 
$\Theta_{1,s}\in \Hom_{inj} (TD,TX|_{G_1(D)})$ $(s\in[0,1])$
satisfying
\[
	\Theta_{1,0}=dG_1\colon TD\to TX|_{G_1(D)}, \quad
	\Theta_{1,1} = \phi_1|_{TD} \colon TD \to TX|_{G_1(D)}. 
\]
The homotopy $\Theta_{1,s}$ is fixed over $S$ where it coincides with 
$\wt \theta_{1,s}= dG_1|_{TD|_S}$ by property (iv'). 
(On Fig.\ \ref{Fig4}, $\Theta_{1,s}$ appears on the right face of the cube, 
with bold vertical sides indicating that it is constant on 
$TD|_S$ where it equals $\ell_1\oplus dg_1$.)

%
%
%
%

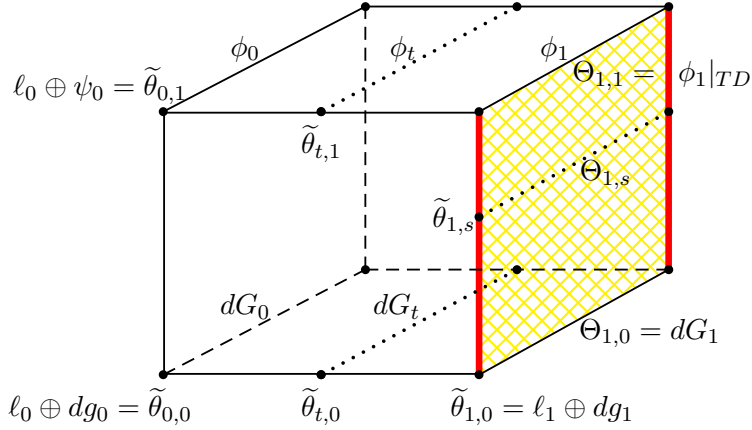
\begin{figure}[ht]
\psset{unit=0.7cm, xunit=1.2, linewidth=0.7pt}  

\begin{pspicture}(-1,-1)(9,8)

\pspolygon[linestyle=none,fillstyle=crosshatch,hatchcolor=yellow](5,0)(5,5)(8,7)(8,2)  

\psframe(0,0)(5,5)                                                       

\psline(5,0)(8,2)                                                        
\psline(5,5)(8,7)
\psline(8,2)(8,7)
\psline[linewidth=2.5pt,linecolor=red](5,0)(5,5)
\psline[linewidth=2.5pt,linecolor=red](8,2)(8,7)

\psline(0,5)(3.2,7)                                                      
\psline(3.2,7)(8,7)

\psline[linestyle=dashed](0,0)(3.2,2)                                    
\psline[linestyle=dashed](3.2,2)(8,2)
\psline[linestyle=dashed](3.2,2)(3.2,7)

\psline[linestyle=dotted, linewidth=1.4pt](2.5,0)(5.6,2)                 
\psline[linestyle=dotted, linewidth=1.4pt](2.5,5)(5.6,7)
\psline[linestyle=dotted, linewidth=1.4pt](5,3)(8,5)

\psdots(0,0)(2.5,0)(5,0)(0,5)(2.5,5)(5,5)                                
(3.2,2)(8,2)(3.2,7)(8,5)(2.5,0)(5.6,2)(2.5,5)(5.6,7)(8,7)(5,3)

\rput(-1,-0.6){$\ell_0\oplus dg_0=\wt\theta_{0,0}$}                      
\rput(2.5,-0.6){$\wt\theta_{t,0}$} 
\rput(6,-0.6){$\wt\theta_{1,0}=\ell_1\oplus dg_1$} 
\rput(-1,5.5){$\ell_0\oplus\psi_0=\wt\theta_{0,1}$}
\rput(2.5,4.4){$\wt\theta_{t,1}$}

\rput(1.3,1.3){$dG_0$}
\rput(3.7,1.3){$dG_t$}

\rput(1.3,6.2){$\phi_0$}
\rput(3.8,6.2){$\phi_t$}
\rput(6.2,6.2){$\phi_1$}
\rput(7.9,5.7)  {$\Theta_{1,1}=\ \,\phi_1|_{TD}$}            

\rput(7,3.8){$\Theta_{1,s}$}
\rput(7.7,0.8){$\Theta_{1,0}=dG_1$}
\rput(4.6,3){$\wt\theta_{1,s}$}

\end{pspicture}
\caption{The homotopy $\Theta_{t,s}$}
\label{Fig4}
\end{figure}

Since $\phi_1 \colon T\C^n|_D \to TX|_{G_1(D)}$ 
is a $\C$-linear vector bundle isomorphism, the h-principle for 
totally real immersions (see \cite{EM,Ghprinciple,Gbook})
provides a regular homotopy of immersions $G_t\colon D\to X\bs W$ ($t\in [1,2]$), 
fixed near $S$, such that the immersion $G_2$ is $J$-real and its 
complexified differential $d^\C G_2$ is homotopic to $\phi_1$ in the space
of $\C$-linear maps $T\C^n|_D \to TX$ of maximal rank.
If in addition $G_1$ is an embedding,
we can deform it to a totally real embedding $G_2$
by an isotopy which is fixed near $S$;
this follows from the fact that totally real embeddings 
also satisfy the h-principle (see \cite{Gbook}). 
For $k<n$ or $k=n>2$ this can also be seen 
from the results in \cite{Ftotallyreal}, and for $k=n=2$ 
it follows from the work of Eliashberg and Harlamov \cite{EH} on 
cancellation of complex points of real surfaces in complex surfaces
(this is discussed in \S 7 below; see also \cite{Fcomplexpoints}).

Finally we reparametrize the family 
$\{G_t \colon t\in [0,2]\}$ back to the parameter interval $[0,1]$.
This proves the existence of a regular homotopy with the required properties.

It remains to be seen that, unless $k=n=2$, there also exists an 
{\em isotopy of embeddings}  $\{G_t\}$ with these properties. 
If $k<n$, a small perturbation of $\{g_t\}$ 
with fixed ends at $t=0,1$ gives an isotopy which can be realized by 
an ambient diffeotopy, and we get an isotopy of embedded discs 
$G_t\colon D\hra X\bs W$ with $G_t|_S=g_t$. 
For $k=n=1$ the conclusion of Lemma \ref{L3.1} 
obviously holds for any attached 1-disc (segment).

%
%
%
%
In the remainder of the proof we consider the case $k=n>2$.
The main idea of the following argument was communicated to us 
by Y.\ Eliashberg (personal communication, June 2005).

A generic choice of the regular homotopy $g_t\colon S\simeq S^{n-1}\to \Sigma$
insures that $g_t$ is an embedding for all but finitely many
parameter values $t\in [0,1]$, and it has a simple (transverse) double point 
at each of the exceptional parameter values. We wish  to change 
the Legendrian embedding $g_1$ by a regular homotopy
of Legendrian immersions $g_t\colon S\hra \Sigma$ $(t\in[1,2])$ 
to another Legendrian embedding $g_2$ so that the 
resulting regular homotopy $\{g_t\colon t\in [0,2]\}$ will have 
self-intersection index zero. More precisely, the map
$\wt g\colon \wt S=S\times [0,2] \to \wt \Sigma= \Sigma\times[0,2]$,
defined by $\wt g(x,t)=(g_t(x),t)$, is an immersion 
of the $n$-dimensional oriented manifold $\wt S$ into the $2n$-dimensional
oriented manifold $\wt \Sigma$ such that the double points of $\wt g$ 
correspond to the double points of the regular homotopy $\{g_t\}$,
and we define the index $i(\{g_t\})$ as the number 
of double points of $\wt g$, counted with their orientation signs. 

If this index equals zero then a foliated version of 
the Whitney trick allows us to deform $\{g_t\}_{t\in[0,2]}$
with fixed ends to an isotopy of embeddings. 
This is done by connecting a chosen pair of double points
$q_0,q_1\in \wt g(\wt S)$ of the opposite sign, 
lying over two different values $t_0<t_1$ of the parameter, 
by a pair of curves $\lambda_j(t) = \wt g(c_j(t),t)$
($t\in [t_0,t_1]$, $j=1,2$) which together bound an embedded 
Whitney disc $D^2 \subset \wt \Sigma$ such that $D^2\cap (\Sigma_t\times\{t\})$ 
is an arc connecting $\lambda_1(t)$ to $\lambda_2(t)$ 
for every $t\in [t_0,t_1]$, and it degenerates to 
$q_0$ resp.\ $q_1$ over the endpoints $t_0$ resp.\ $t_1$. 
The rest of the procedure, removing this pair of double
points by pulling $\wt g(\wt S)$ across $D^2$, is standard \cite{Whitney}.
Performing this operation finitely many times one can remove all
double points and change $\{g_t\}$ to an isotopy of embeddings. 

The rest of the proof can be completed
exactly as before: we extend $g_t$ to an isotopy of embedded
discs $G_t\colon D\hra X\bs W$, with $G_t|_S=g_t$,
covered by a homotopy of $\C$-linear isomorhisms 
$\phi_t \colon T\C^n|_D\to TX|_{G_t(D)}$.
Observe that $\{dg_t\}$ still has the correct homotopy
property so that the final embedding $G_2$ can be deformed
(with fixed boundary) to a totally real embedding.  

It remains to see that the index of $\{g_t\}_{t\in[0,1]}$ 
can be changed to an arbitrary number (in particular, to zero) 
by a small Legendrian deformation of $g_1$ in $\Sigma$.
This will be realized by a local Legendrian isotopy
which introduces the correct number of double points.
(A similar deformation is used in \cite[\S 2.4]{E} for 
changing a stably special HAT to a special one.)

Set $L=L_1:=g_1(S)\subset \Sigma$, an embedded Legendrian sphere. 
Choose a point $a\in L\subset \Sigma$. In suitable local 
coordinates $(z,q,p)\in \R^{2n-1}$ on $\Sigma$, with $a$ 
corresponding to $0\in\R^{2n-1}$, the contact form is 
$\eta=dz - \sum_{j=1}^{n-1} p_j dq_j$, and $L$ is given by 
the equations
\[
	\{(z,q,p)\in\R^{2n-1} \colon 
	z^2=q_1^3,\ p_1^2= \frac{9}{4}q_1,\ p_2=\cdots=p_{n-1}=0\}.
\]
(Section 2.4 in \cite{E}.)
Let $\pi \colon\R^{2n-1}\to\R^{n-1}$ denote the projection $\pi(z,q,p) =  q$. 
Choose a closed ball $\Delta \subset \R^{n-1}$ centered at  
$(q_1^0,0,\ldots,0)$ for a small $q_1^0>0$,  of radius $q_1^0/2$.  
Let $\phi \colon \Delta\to \R$ be a smooth function 
equaling $0$ near $\di\Delta$. Set
\[
	h_t(q) =q_1^{3/2}\bigl(1 + (t-1) \phi(q)\bigr), \quad t\in [1,2].
\]
Let $L_t$ equal $L$ outside of $\pi^{-1}(\Delta)$ and equal  
\[
	\bigl\{(z,q,p) \colon z = h_t(q),\ p=\frac{\di h_t}{\di q}(q) \bigr\}
	\cup \bigl\{(z,q,p) \colon z=-h_t(q),\ p=-\frac{\di h_t}{\di q}(q) \bigr\}
\]
over $\Delta$. (We choose $\phi$ with sufficiently small
derivative to insure that we remain in the given coordinate patch; 
this can be done if $q_1^0>0$ is chosen small enough.)
Let $g_t\colon S\to \Sigma$ $(t\in [1,2])$ be the regular homotopy 
such that $g_t(S)=L_t$. 
The deformation is illustrated by Fig.\ \ref{vozel}.
The top diagrams show the projection onto the $(z,q)$-plane 
at three typical stages, with the cusp at $(z,q)=(0,0)$  
and with a self-intersection shown in the middle figure.

%
%
%
%
\begin{figure}[ht]
\psset{unit=0.6cm, linewidth=0.7pt} 
\begin{pspicture}(-8,-4)(8,4)

\psecurve
(-100,2)(-8,2)(-7,2.125)(-6,2.35)(-5,2.65)(-4,3)(-4,3)
\psecurve
(-100,2)(-8,2)(-7,1.875)(-6,1.65)(-5,1.35)(-4,1)(-4,1)

\psecurve(-100,2)(-2,2)(-1,2.125)(0,2.35)
\psecurve(-2,2)(-1,2.125)(-0.5,2)(0,2.35)(1,2.65)
\psecurve(-1,2.125)(0,2.35)(1,2.65)(2,3)(2,3)

\psecurve(-100,2)(-2,2)(-1,1.875)(0,1.65)
\psecurve(-2,2)(-1,1.875)(-0.5,2)(0,1.65)(1,1.35)
\psecurve(-1,1.875)(0,1.65)(1,1.35)(2,1)(2,1)

\psecurve(-100,2)(4,2)(5,2.125)(6,2.35)
\psecurve(4,2)(5,2.125)(5.5,1.9)(6,2.35)
\psecurve(5,2.125)(6,2.35)(7,2.65)(8,3)(8,3)

\psecurve(-100,2)(4,2)(5,1.875)(6,1.65)
\psecurve[doubleline=true,doublecolor=black, linewidth=1.5pt,
                   linecolor=white,                     
                   doublesep=0.7pt](4,2)(5,1.875)(5.5,2.1)(6,1.65)
\psecurve(5,1.875)(5.5,2.1)(6,1.65)(7,1.35)
\psecurve(5,1.875)(6,1.65)(7,1.35)(8,1)(8,1)

\psecurve[doubleline=true,doublecolor=black, linewidth=1.5pt,
                   linecolor=white,                     
                   doublesep=0.7pt](5,2.125)(5.5,1.9)(6,2.35)(7,2.65) 

\pscurve[linewidth=0.2pt]{->}(-3.8,2.5)(-3,2.7)(-2.2,2.5)
\pscurve[linewidth=0.2pt]{->}(2.2,2.5)(3,2.7)(3.8,2.5)

\pscurve[doubleline=true,doublecolor=black, linewidth=2pt,
                   linecolor=white,                     
                   doublesep=2pt](-4.4,-1.2)(-8,-2)(-4,-3)

\psecurve[doubleline=true,doublecolor=black, linewidth=2pt,
                   linecolor=white,                     
                   doublesep=2pt](-0.5,-2.2)(-1.1,-1.6)(-2,-2)(-1,-2.6)(-0.5,-1.8)

\psecurve[doubleline=true,doublecolor=black, linewidth=2pt,
                   linecolor=white,                     
                   doublesep=2pt](-2,-2)(-1,-2.6)(-0.5,-1.8)(0,-2.8)

\psecurve[doubleline=true,doublecolor=black, linewidth=0pt,
                   linecolor=white,                     
                   doublesep=2pt](-1,-2.6)(-0.5,-1.8)(0,-2.8)(2,-3)

\psecurve[doubleline=true,doublecolor=black, linewidth=2pt,
                   linecolor=white,                     
                   doublesep=2pt](-2,-2)(-1.1,-1.6)(-0.5,-2.2)(0,-1.5)

\psecurve[doubleline=true,doublecolor=black, linewidth=0pt,
                   linecolor=white,                     
                   doublesep=2pt](-1.1,-1.6)(-0.5,-2.2)(0,-1.5)(1.7,-1.2)

\psecurve[doubleline=true,doublecolor=black, linewidth=2pt,
                   linecolor=white,                     
                   doublesep=2pt](3,-1.2)(1.7,-1.2)(0,-1.5)(-0.5,-2.2)

\psecurve[doubleline=true,doublecolor=black, linewidth=2pt,
                   linecolor=white,                     
                   doublesep=2pt](3,-3)(2,-3)(0,-2.8)(-0.5,-1.8)

\psecurve[doubleline=true,doublecolor=black, linewidth=2pt,
                   linecolor=white,                     
                   doublesep=2pt](5.5,-2.2)(4.9,-1.6)(4,-2)(5,-2.6)(5.5,-1.7)

\psecurve[doubleline=true,doublecolor=black, linewidth=2pt,
                   linecolor=white,                     
                   doublesep=2pt](4,-2)(5,-2.6)(5.5,-1.7)(6,-2.8)

\psecurve[doubleline=true,doublecolor=black, linewidth=2pt,
                   linecolor=white,                     
                   doublesep=2pt](4,-2)(4.9,-1.6)(5.5,-2.4)(6,-1.5)

\psecurve[doubleline=true,doublecolor=black, linewidth=2pt,
                   linecolor=white,                     
                   doublesep=2pt](4.9,-1.6)(5.5,-2.4)(6,-1.5)(7.7,-1.2)

\psecurve[doubleline=true,doublecolor=black, linewidth=2pt,
                   linecolor=white,                     
                   doublesep=2pt](5,-2.6)(5.5,-1.7)(6,-2.8)(8,-3)

\psecurve[doubleline=true,doublecolor=black, linewidth=2pt,
                   linecolor=white,                     
                   doublesep=2pt](9,-1.2)(7.7,-1.2)(6,-1.5)(5.5,-2.4)

\psecurve[doubleline=true,doublecolor=black, linewidth=2pt,
                   linecolor=white,                     
                   doublesep=2pt](9,-3)(8,-3)(6,-2.8)(5.5,-1.7)

\pscurve[linewidth=0.2pt]{->}(-3.8,-1.5)(-3,-1.3)(-2.2,-1.5)
\pscurve[linewidth=0.2pt]{->}(2.2,-1.5)(3,-1.3)(3.8,-1.5)

\end{pspicture}

\caption{Changing the index of a regular homotopy by $+1$.}
\label{vozel}

\end{figure}
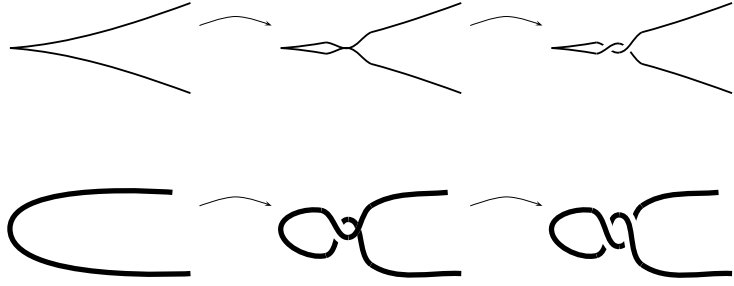

The index of $\{g_t\}_{t\in [1,2]}$ equals the intersection number 
between the manifolds 
\[
	M_\pm =\bigl\{(z,q,p,t)\in\R^{2n} \colon 
	z= \pm h_t(q),\ p=\pm \frac{\di h_t}{\di q}(q),
	\ q\in\Delta,\  t\in [1,2] \bigr\}.
\]
The intersection points of $M_+$ and $M_-$ are solutions of the 
equations  
\[
	1+(t-1)\phi=0,\quad \frac{\di\phi}{\di q}=0,\qquad t\in [1,2]. 
\] 
This is precisely the set of all critical points of $\phi$ with
the critical values belonging to $(-\infty,-1]$.
By a generic choice of $\phi$ we can insure 
that $-1$ is not a critical value of $\phi$. 
A computation shows that each point $(q,t)$ satisfying the 
above equations adds $\pm 1$ to the index $i(\{g_t\})$, 
depending on the sign of the determinant of 
the Hessian $\mathrm Hess(\phi)$ at $q$; hence we get $+1$ 
at a critical point  of even Morse index and $-1$ at a
critical point  of odd Morse index. 
Similarly, as we increase $c\in\R$, the Euler characteristic 
of the sublevel set $\{\phi \le c\}$
increases by one at every critical point of $\phi$ of even Morse index,
and it decreases by one at every critical point of odd Morse index.
We conclude that $i=i(\{g_t\}_{t\in [1,2]})$ equals
the Euler number of the set $\{q\in \Delta \colon \phi(q)\le -1\}$. 
If $n>2$, this can be arranged to equal any preassigned integer
by a suitable choice of $\phi$, and hence we can arrange the index
$i$ to equal zero. If $n=2$ then $i$ can be arranged to be any nonnegative number 
since  $\{\phi\le -1\}$ is a union of segments, 
but it cannot be negative.

%
This completes the proof of Lemma \ref{L3.1} in the smooth case. 
Assume now that $J$ is integrable in a neighborhood of the disc 
$G_0(D)$ and the hypersurface $\Sigma=\di W$ is real analytic near the attaching 
sphere $G_0(S)\subset \Sigma$ with respect to the induced complex structure;
we wish to find a real analytic disc satisfying the conclusion of Lemma \ref{L3.1}.
Since the disc $G_1$ constructed above can be chosen arbitrarily 
$\cC^0$-close to $G_0$, we may assume that the same conditions 
on $J$ and $\Sigma$ also hold near $G_1(D)$.
By \cite[Lemma 2.5.1.]{E} (which uses Gray's theorem on real analytic 
approximation  of Legendrian embeddings) 
it is possible to approximate $G_1$ in the $\cC^1$ topology
by a disc $G'_1\colon (D,S)\to (X\bs W,\Sigma)$ which is real analytic
near $S$ such that $G'_1|_S\colon S\hra \Sigma$ is Legendrian
and $G'_1$ is normal to $\Sigma$ along $S$. 
It remains to perturb $G'_1$ to a nearby real analytic map 
$G''_1\colon D\to X$ which agrees with $G'_1$ to the second order along $S$
and to replace $G_1$ by $G''_1$.  
\end{proof}

%
%
%
%
\begin{remark}
\label{framing}
We wish to point out that in the critical case $k=n>2$ 
Lemma \ref{L3.1} does not seem a direct consequence of the results in \cite{E}.
Let $D=D^n$ and $S=\di D$. In \cite[\S 2]{E} it is shown that for $n\ne 2$ 
every embedding $G\colon (D,S) \to (X\bs W,\Sigma)$ can be isotoped to one 
for which $G$ is totally real near the boundary, 
the attaching sphere $G(S)$ is Legendrian in $\Sigma$, 
and the associated normal framing $\beta\colon\nu\to TX|_{G(D)}$ 
satisfies 
\begin{equation}
\label{beta}
	\beta \circ J_{st} = J\circ dG \ \ \text{on}\ T_x D, \quad x\in S. 
\end{equation}
This means that the triple $(G,\beta,\phi=dG\oplus \beta)$, restricted to the
points of $S$, is a {\em special HAT} in the sense of \cite[p.\ 33]{E}. 
Choose a diffeomorphism $\wt G$ from a standard handle 
$H\subset\C^n$  onto $\wt G(H)\subset X$ such that $\wt G|_D=G$ 
and $d\wt G|_\nu =\beta$. By (\ref{beta}) the 
push-forward $\wt J = \wt G_*(J_{st})$ of the standard complex structure 
on $\C^n$ agrees with $J$ at every point of $G(S)$, and by a small correction 
these structures can be made to agree in a neighborhood of $G(S)$ in $X$.
Using the h-principle and some homotopy theory one can see 
that the disc $G$ (which is clearly $\wt J$-real)
is isotopic to a $J$-real disc by an isotopy which is fixed near $S$ 
precisely when $\wt J$ is homotopic to $J$ along $G(D)$ by a homotopy 
which is fixed near $G(S)$. The latter condition holds if and only if 
the topological invariant in (\ref{difference}) below vanishes. 
The following proposition shows that this always holds
for certain values of $n$ modulo 8, and the proof of 
Lemma \ref{L3.1} avoids this potential 
problem in every dimension.
\end{remark}

\begin{proposition}
\label{P3.4}
Let $W$ be an open, relatively compact domain with smooth 
strongly pseudoconvex boundary $\Sigma=\di W$ in an 
almost complex manifold $(X,J)$ of real dimension $2n$.
Let $D=D^n$, $S=\di D$, and let $G\colon (D,S)\to (X\bs W,\Sigma)$
be a smooth embedding, normal to $\Sigma$, such that $G(S)$ is 
Legendrian in $\Sigma$. Assume that $G$ admits a normal framing $\beta$ 
satisfying (\ref{beta}) over $S$. If $n\in \{1,3,4,5\}$ modulo 8 
then $G$ is isotopic to a totally real embedding by an isotopy 
which is fixed near $S$.
\end{proposition}

\begin{proof} 
We extend $G$ to a diffeomorphism $H\to G(H)\subset X$ satisfying $dG|_\nu=\beta$. 
The difference on $D$ between the almost complex structures 
$J' = \wt G^*J$ and $J_{st}$ (which agree over $S$) defines an element
\begin{equation}
\label{difference}
	\delta(J',J_{st}) \in [S^n,GL^+_{2n}(\R)/GL_n(\C)] = \pi_n(SO(2n)/U(n))
\end{equation}
(Milnor \cite[p.\ 133]{Mi}). Using the long exact sequence of 
homotopy groups and the five lemma one sees that 
\[
		\pi_n(SO(2n)/U(n)) = \pi_n(SO/U)=\pi_n(\Omega SO)=\pi_{n+1}(SO).
\]
By the (real) Bott periodicity theorem this group equals 
$\Z$ if $n\in \{2,6\}$ modulo 8, it equals $\Z_2$ if $n\in \{0,7\}$ modulo 8, 
and it vanishes for the remaining values $n\in\{1,3,4,5\}$ modulo 8.
In the last case we conclude that $J=J_1$ is homotopic along $G(D)$ to 
$J_0= G_*(J_{st})$ by a homotopy which is fixed near 
$G(S) \subset \Sigma$. Since $G_0=G$ is clearly $J_0$-real, 
Gromov's h-principle for totally real embeddings \cite{EM,Ftotallyreal,Gbook}
gives an isotopy of embedded discs $G_t\colon (D,S)\to (X\bs W,\Sigma)$
which is fixed near $S$ such that $G_t$ is $J_t$-real for every $t\in[0,1]$.
At $t=1$ we get a $J$-real embedded disc~$G_1$.
\end{proof}

%
%
%
%
\section{A holomorphic approximation theorem}
Let $X$ be a complex manifold with an integrable
complex structure $J$. We denote by 
$\cH(X)=\cH(X,J)$ the algebra of all holomorphic 
functions on $X$. A compact set $K$ in $X$ is $\cH(X)$-convex  
(or $\cH(X,J)$-convex) if for every point $p\in X\bs K$ there exists 
an $f\in \cH(X)$ with $|f(p)|> \sup_{x\in K} |f(x)|$.

We say that a compact set $K$ in $X$
is {\em holomorphically convex} if there is an 
open Stein domain $\Omega\subset X$ 
containing $K$ such that $K$ is $\cH(\Omega)$-convex.
By the classical theory \cite[Chapter 2]{Ho} this is equivalent 
to the existence of a Stein neighborhood $\Omega$ of $K$ and a continuous
plurisubharmonic function $\rho\ge 0$ on 
$\Omega$ such that $\rho^{-1}(0)=K$ and  
$\rho$ is strongly plurisubharmonic on $\Omega\bs K$. 
We may take $\Omega = \{\rho<c_1\}$ for some $c_1>0$; 
for any $c\in (0, c_1)$ the sublevel set 
$\{\rho<c\} \Subset \Omega$ 
is then Stein and Runge in $\Omega$ \cite[\S 4.3]{Ho}.

A point $p_0$ in an immersed real $k$-dimensional
submanifold $M \subset X$ is said to be a {\em special double point}
if there is a holomorphic coordinate system
\[
	z=(x'+iy',z'') \colon U\to \wt U\subset \C^n = \R^k\oplus i\R^k\oplus \C^{n-k}
\]
in a neighborhood $U\subset X$ of $p_0$ such that $z(p_0)=0$ and
\begin{equation}
\label{specialdouble}
	z(M\cap U) = \wt U\cap \left ( \{(x'+i0',0'')\colon x'\in\R^k\} \cup
	              \{(0'+iy',0'')\colon y'\in \R^k\} \right)  
\end{equation}	

The following approximation  theorem will be used in the proof of 
Lemma \ref{L5.1}. It is far from the most general one with respect 
to the type of allowed double points of $M$, but it will suffice 
for our purposes.

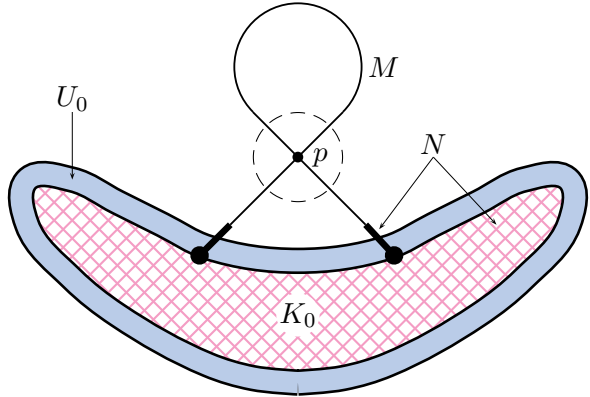
\begin{figure}[ht]
\psset{unit=0.6cm, linewidth=0.7pt} 

\begin{pspicture}(0,0)(16,10)

\pscurve[fillstyle=crosshatch,hatchcolor=Lavender,linewidth=1pt,
doubleline=true,doublesep=8pt,dimen=outer,doublecolor=CornflowerBlue]
(8,1)(9,1.1)(10,1.4)(12,2.5)(13,3.3)(14,4.5)(14,5.5)(13,5.5)(12,5)(11,4.5)(10,4)(8,3.7) 
(6,4)(5,4.5)(4,5)(3,5.5)(2,5.5)(2,4.5)(3,3.3)(4,2.5)(6,1.4)(7,1.1)(8,1)

\psdots[dotscale=1.2](8,6)

\psline(6.5,4.5)(9,7)
\psline(7,7)(9.5,4.5)
\psarc(8,8){1.41}{-45}{225}
\pscircle[linewidth=0.4pt,linestyle=dashed](8,6){1}         

\psline[linewidth=2.5pt,linecolor=black]{*-}(5.82,3.82)(6.5,4.5)    
\psline[linewidth=2.5pt,linecolor=black]{-*}(9.5,4.5)(10.13,3.82)   

\psline[linewidth=0.2pt]{<-}(9.8,4.4)(11,6)                 
\psline[linewidth=0.2pt]{->}(11,6)(12.4,4.5)                
\rput(11,6.3){$N$} 					    

\psline[linewidth=0.2pt]{->}(3,7)(3,5.5)
\rput(3,7.3){$U_0$} 

\rput(8.5,6){$p$}
\rput(9.9,8){$M$}

\pscircle[fillstyle=solid,fillcolor=white,linestyle=none](8,2.5){0.5}
\rput(8,2.5){$K_0$}

\end{pspicture}

\caption{A kinky disc $M$ attached to $K_0$}
\label{Fig5}
\end{figure}

\begin{theorem} 
\label{T4.1}
Let $K_0$ and $M$ be compact sets in a complex manifold $X$, 
where $M\bs K_0$ is a smoothly immersed totally real submanifold 
such that each non-smooth point is a special double point (\ref{specialdouble}). 
Assume that $K_0$ is holomorphically convex and 
there is a compact holomorphically convex relative neighborhood $N$ 
of $K_0$ in $K:=K_0\cup M$ (Fig.\ \ref{Fig5}).

Given a continuous map $f\colon X\to Y$ to a complex manifold $Y$ 
such that $f$ is holomorphic in an open neighborhood of $K_0$, there exist 
open Stein domains $V_1\supset V_2\supset \cdots \supset \cap_j V_j=K$ 
and holomorphic maps $f_j\colon V_j\to Y$ $(j=1,2,\ldots)$ such that
$K_0$ and $K$ are $\cH(V_j)$-convex for every $j$, and 
$f_j|_K \to f|_K$ uniformly as $j\to\infty$. 
If $M\bs K_0$ is an embedded submanifold and 
$f|_{M\bs K_0}$ is of class $\cC^r$ then $f_j$ 
can be chosen such that $f_j|_K\to f|_K$ uniformly 
and $f_j|_{M\bs K_0} \to f|_{M\bs K_0}$ in the $\cC^r(M)$ topology
as $j\to\infty$. 

The analogous results hold for a family of maps 
$f\colon X\times P\to Y$ paramet\-rized by a compact
Hausdorff space $P$.
\end{theorem}

In the case when $M\bs K_0$  is an embedded totally real submanifold 
this is essentially Theorem 3.2 in \cite{FFourier1}.
Uniform approximation was obtained by H\"ormander and Wermer \cite{HW}
in the case $X=\C^n$, $Y=\C$ and $M\bs K_0$ a totally real embedded $\cC^1$ submanifold.
Another special case is due to Forn\ae ss and Stout \cite{FS}.

%
%
%
%
%
%

\begin{proof}
Consider first the case when $M\bs K_0$ is an {\em embedded} totally real 
submanifold. The assumption regarding $N$ 
implies by \cite[Theorem 3.1]{FFourier1} that the set 
$K=K_0\cup M$ is holomorphically convex. (See also
\cite{HW} for the case $X=\C^n$.)
More precisely, given an open neighborhood $U_0\subset X$ of $K_0$,
there exists a continuous plurisubharmonic function $\rho\ge 0$ 
in an open neighborhood $U\subset X$ of $K$ such that 
$\rho$ is smooth strongly plurisubharmonic on $U\bs U_0$ 
and $\rho^{-1}(0)=K$ \cite[p.\ 1923]{FFourier1}. 
Taking $V_j=\{x\in U\colon \rho(x)<c_j\}$ for a decreasing 
sequence of small numbers $c_j>0$ converging to zero we get 
a Stein neighborhood basis of $K$ such that
$K$ is $\cH(V_j)$-convex for each $j$.
It is easily seen that $K_0$ is also $\cH(V_j)$-convex. 
Indeed, choose a smooth function $\chi\ge 0$ 
with compact support in $U$ which vanishes on an open set $U_1\supset \overline U_0$ 
and is positive on $M\bs \overline U_1$. Then the function $\rho+\epsilon \chi\ge 0$
is still plurisubharmonic for sufficiently small $\epsilon>0$
and hence its zero set $K\cap \overline U_1$ is $\cH(K)$-convex
(i.e., for every point $p\in K\bs \overline U_1$ there is a holomorphic function
$g$ in a neighborhood of $K$ with $1=g(p)>\sup_{x\in K\cap \overline U_1} |g(x)|$).
Since $U_0$ and $U_1$ can be chosen arbitrary close to $K_0$,
we see that $K_0$ is $\cH(K)$-convex, and hence also $\cH(V_j)$-convex.

By approximation we may assume that $f$ is smooth on $M\bs K_0$.
Theorem 3.2 in \cite{FFourier1} now shows 
that $f$ can be approximated uniformly on $K_0$, an in
the $\cC^r$ topology on $M\bs K_0$, by maps which are holomorphic
in small open neighborhoods of $K$ in $X$
(the size of the neighborhood depends on the rate of approximation).
Choosing these neighborhoods from the above sequence $V_j$ 
we get the conclusion of Theorem \ref{T4.1}.

Suppose now that $M\bs K_0$ is an {\em immersed} totally
real submanifold with special double points
$p_1,\ldots,p_m$. (Replacing $K_0$ by a relative 
neighborhood of $K_0$ in $K$ we can 
assume that there are only finitely many such points.) 
Let $B_j \subset X\bs K_0$ be a small open neighborhood
of $p_j$ in $X$ such that $\overline B_j$ is mapped onto a closed ball 
around $0 \in \C^n$ by a local coordinate map (\ref{specialdouble}).
By a uniformly small change we make $f$ smooth on $M\bs K_0$ 
and constantly equal to $f(p_j)$ on a neighborhood
of $\overline B_j$ in $X$; the latter change 
can be made small by choosing the balls $B_j$ as small as necessary. 
We can now apply the previous argument with the set
$K'_0=K_0\cup (\cup_{j=1}^m \overline B_j)$
and the embedded totally real submanifold $M'=M\bs K'_0$. 
Indeed,  $K'_0$ is clearly holomorphically convex,
and it has a compact holomorphically convex relative neighborhood $N'$ in 
$K':= K'_0\cup M'= K\cup (\cup_{j=1}^m \overline B_j)$.
(It suffices to take $N'=N\cup(\cup_{j=1}^m  N_j)$ where 
$N_j$ is the union of $\overline B_j$ with a suitably chosen 
small collar in $M$; in local coordinates, $N_j$ 
corresponds to the union of a closed ball in $\C^n$ centered at $0$ 
with a collar in $\R^k\cup i\R^k\subset\C^k\times\{0\}^{n-k}$.)
As before, Theorem 3.1 in \cite{FFourier1} implies that
$K'$ is holomorphically convex, and Theorem 3.2 in \cite{FFourier1} 
gives a desired approximation of $f$ by holomorphic maps in small Stein 
neighborhoods of $K'$ in $X$.

The proof of Theorem \ref{T4.1} for a family of maps 
$f\colon X\times P\to Y$, with $P$ a compact Hausdorff space, 
is obtained by covering the graph of the family in $X\times Y$ 
(after an initial smoothing of the maps $f(\cdotp,p) \colon X\to Y$
on the $M\bs K_0$) by finitely many 
Stein neighborhoods in $X\times Y$, 
using these to approximate $f$ by local (in $P$)
families of holomorphic maps, and patching these families
by a continuous partition of unity in the parameter $p\in P$.
The latter is possible since we can introduce a complex 
linear structure on the fibers of the projection $X\times Y\to X$ 
within a small Stein neighborhood of each individual graph. 
The details in a very similar context can be found 
in \cite{FP1} (proof of Theorem 4.2, pp.\ 138-139).
\end{proof}

%
%
%
%
\section{Extending a holomorphic map across a handle}
The following lemma is the key ingredient in the proofs of our main results. 

%
%
%
%
\begin{lemma}
\label{L5.1}
Let $(X,J)$ be an almost complex manifold of real dimension 
$2n$. Let $W\Subset X$ be a smoothly bounded 
domain such that $J$ is integrable in a neighborhood 
of $\overline W$, the manifold $(W,J)$ is Stein, 
and $\Sigma=\di W$ is strongly $J$-pseudoconvex. 
Let $D=D^k$ and $S=S^{k-1}=\di D$ $(1\le k \le n)$. 
Let $G\colon (D,S) \to (X\bs W,\Sigma)$ 
be a smooth $J$-real embedding which is normal to $\Sigma$ 
and such that $G|_S\colon S\to \Sigma$ is Legendrian.
Assume that $Y$ is a complex manifold and 
$f\colon X\to Y$ is a continuous map which is 
$J$-holomorphic in an open neighborhood of $\overline W$.
Let $d_Y$ be a distance function on $Y$ induced by a 
smooth Riemannian metric. 

After a small smooth perturbation of $G$ there exist  an integrable 
complex structure $\widetilde J$ in an open neighborhood 
$U\subset X$ of $K :=\overline W \cup G(D)$, a homotopy 
$J_t$ $(t\in [0,1])$ of almost complex structures on $X$
which is fixed on a neighborhood of $\overline W$ and on $M=G(D)$,
with $J_0=J$ and $J_1=\wt J$, and for every $\e>0$ there exist 
a smoothly bounded strongly $\widetilde J$-pseudoconvex 
Stein domain $\widetilde W$ and a map $\wt f\colon X\to Y$ 
satisfying the following:
\begin{itemize}
\item[(i)]  $\overline W \cup G(D) \subset \widetilde W \subset U$, 
$\wt W$ is a handlebody with core $K=\overline W \cup G(D)$, and 
$\overline W$ is $\cH({\widetilde W},\widetilde J)$-convex. 
\item[(ii)] The map $\widetilde f|_{\wt W} \colon \wt W\to Y$ 
is $\widetilde J$-holomorphic.
\item[(iii)] 
There is a homotopy $f_t\colon X\to Y$ $(t\in [0,1])$, with  
$f_0=f$ and $f_1=\widetilde f$, such that for each 
$t \in [0,1]$ the map $f_t$ is $J$-holomorphic on 
a neighborhood of $\overline W$ and 
$\sup_{x\in \overline W} d(f(x),f_t(x)) <\e$. 
\end{itemize}

If in addition $f$ is covered by a complex vector bundle map 
$\iota\colon (TX,J)\to TY$ which is of maximal rank on every fiber and
such that $df=\iota$ on a neighborhood of $\overline W$
then we can choose $\wt f$ to be of maximal rank at every point 
of $\wt W$ and such that $d\wt f$ is homotopic to $\iota$ 
through complex vector bundle maps $\iota_t \colon (TX,J_t) \to TY$ 
of pointwise maximal rank.

If $\Sigma$ is real analytic and the almost complex structure 
$J$ is integrable in a neighborhood of $K$ then the above conclusions 
hold with $J=\wt J$. 

The analogous results hold for a continuous family of maps 
with a parameter in a compact Hausdorff space.
\end{lemma}

\begin{proof}
After a small enlargement of $W$ and a small deformation of $G$
we may assume that $\di W$ is real analytic
and strongly $J$-pseudoconvex, $J$ is integrable in a neighborhood 
of $\overline W$, and the $k$-disc $M:=G(D)$ is attached 
to $\overline W$ along the Legendrian $(k-1)$-sphere
$G(S) \subset \di W$. 
By \cite[Lemma 2.5.1.]{E} (which uses Gray's theorem on approximation
of Legendrian embeddings by real analytic Legendrian embeddings) we can approximate 
$G$ by a map which is normal to $\Sigma$ and 
real analytic in a neighborhood of $S=\di D$,  
such that the attaching sphere $G(S)\subset\Sigma$ is 
Legendrian in $\Sigma=\di W$. 

We first consider the case $k=n$. 
For every $x\in D$ let $A_x\colon T_x\C^n \to T_{G(x)} X$
denote the unique $(J_{st},J)$-linear map which agrees with
$dG_x$ on $T_x D$. We extend $G$ to a smooth diffeomorphisms
$\widetilde G$ from a standard handle $H\subset \C^n$ (\ref{Hdelta})
onto a neighborhood $\widetilde H= \widetilde G (H)$ 
of $G(D)$ in $X$ such that $d\wt G_x = A_x$ for each $x\in D$.
Near the sphere $S=\di D$ we take $\wt G$ to be the
complexification of $G$, hence biholomorphic.
If $J$ is integrable then we can choose $G$ to be real analytic 
and $\wt G$ to be its complexification (Lemma \ref{L3.1}).

Let $W'$ be a slightly larger domain in $X$ containing $\overline W$. 
Let $\widetilde J$ denote the complex structure on 
$W_1 := W'\cup \widetilde H$ which equals $J$ on $W'$ and equals 
$\widetilde G_*(J_{st})$ on the handle $\wt H$. 
By choosing the sets $W'\supset \overline W$ and $H\supset D$ sufficiently small
we insure that these two complex structures coincide on 
$W'\cap \widetilde H$ (since $\widetilde G$ maps  
a neighborhood of $S \subset\C^n$ biholomorphically onto
a neighborhood of $G(S)\subset X$).
Notice also that $J=\wt J$ at every point of $M=G(D)$ 
since $d\wt G_x=A_x$ was chosen to be $(J_{st},J)$-linear
for each $x\in D$. This clearly implies the existence of a homotopy 
of almost complex structures $\{J_t\}_{t\in[0,1]}$ on $W_1$ which is fixed
on $W'\cup M$ and satisfies $J_0=J$, $J_1=\wt J$.
If $G$ is real analytic then $J=\wt J$ near $M$ and we can choose 
$J_t=J$ for all $t\in[0,1]$.

Our next goal is to (approximately) extend $f\colon X \to Y$ to a holomorphic
map across the handle. By the assumption $f$ is $J$-holomorphic
in a neighborhood of $\overline W$ in $X$, and we may assume
by approximation that it is smooth on $X$. 
Since $\widetilde J=J$ near $\overline W$, $f$ is also 
$\wt J$-holomorphic near $\overline W$.
We wish to apply Theorem \ref{T4.1} in the complex manifold 
$(W_1,\widetilde J)$, with the compact sets $K_0=\overline W$ 
and $K=K_0\cup M$, in order to obtain a $\widetilde J$-holomorphic map 
$\widetilde f\colon V\to Y$ in an open neighborhood $V\supset K$ 
such that $\widetilde f|_K$ approximates $f|_K$ as close as desired.
In order to do so, we must verify that $\overline W$ has 
a compact holomorphically convex relative neighborhood $N$ in $K$.
It is well known (see e.g.\ \cite[Lemma 1]{Ros}) 
that the problem is local near the attaching sphere $G(S)=\di M\subset \Sigma$.
Thus, taking a closed tubular neighborhood $T\subset X$ of $M=G(D)$,
it suffices to show that the set $\wt G^{-1}(T \cap K) \subset \C^n$
is holomorphically convex for a suitable choice of $T$. 
The latter set is the union of the closed disc $D\subset\R^n \subset \C^n$ 
and a piece of a strongly pseudoconvex domain which essentially looks
like the quadric $Q_\lambda$ (\ref{Qlambda}). In fact, by a small
outward bumping  of $\Sigma=\di W$ (from the side of $W$) which
is localized in a tubular neighborhood of the circle $G(S)$ 
(keeping $\Sigma$ and its tangent bundle fixed on $G(S)$) 
we can reduce to the situation when 
$\wt G^{-1}(T\cap K) = {\wt G}^{-1}(T) \cap (Q_\lambda \cup D)$
and $\wt G^{-1}(T)$ is a compact convex set in $\C^n$. 
The holomorphic polynomial $h(z)=z_1^2+\ldots +z_n^2$ on $\C^n$ maps 
the disc $D$ to the segment $[0,1] \subset \R\subset \C$, it maps
the sphere $S=\di D$ to the point $1$, and $\Re h > 1$ on 
the set $Q_\lambda \bs S$. 
(Compare with \cite{Ros} and the proof of Lemma 6.6 in \cite{FActa}.) 
Thus $h$ separates the polynomially 
convex sets $Q_\lambda \cap \wt G^{-1}(T)$ and $D$, and hence 
their union is polynomially convex in $\C^n$ by a lemma
of Eva Kallin \cite[Lemma 29.21]{St}.

Thus Theorem \ref{T4.1} applies and gives a $\wt J$-holomorphic map 
$\wt f$ in a neighborhood of $K$ which approximates $f$ uniformly on $K$.
A homotopy from $f$ to $\widetilde f$ with the required properties 
clearly exists near $K$ provided that the approximation is  
sufficiently close, and it is then used to patch $\wt f$ with 
$f$ outside of a larger neighborhood of $K$.

Remaining in the case $k=n$ for the moment, we also consider the situation 
when $f$ is of maximal complex rank in a neighborhood of $\overline W$  
and is covered by a complex vector bundle map $\iota\colon (TX,J) \to (TY,J_Y)$ 
of fiberwise maximal rank, with $\iota=df$ near $\overline W$.
In this case we must show that $\wt f$ can also be chosen
of maximal rank on $M$, and hence in a neighborhood of 
$K$ provided that the approximation is sufficiently close.
To this end we first deform $f$ (without changing it near $\overline W$) 
such that for every $x\in M$ its differential $df_x \colon T_x M \to T_{f(x)} Y$ is 
of maximal complex rank (equal to ${\rm min}\{n,\dim Y\}$),
and the map $x\to df_x$ ($x\in M$) is homotopic to 
$x\to \iota|_{T_x M}$ by a homotopy of vector bundle maps 
of pointwise maximal complex rank which is fixed near $\di M$. 
This is a straightforward application of Gromov's h-principle, the main 
point being that the pertinent differential relation is {\em ample} on any 
totally real submanifold. For details see \cite[Lemma 6.4]{FActa}
or \cite[Lemma 4.3, p.\ 1931]{FFourier1}.
Applying Theorem  \ref{T4.1} to this new map $f$ we obtain a
$\wt J$-holomorphic map $\wt f$ in a neighborhood of $K$
which approximates $f$ uniformly on $\overline  W$ and 
in the $\cC^1$-topology on $M$. If the approximation
is sufficiently close, the latter property insures that $\wt f$ is of maximal 
rank at every point of $K$. The existence of a homotopy from 
$\iota$ to $d\wt f$ with the required properties follows 
from the construction (see \cite{FActa} for the details).

To complete the proof of Lemma \ref{L5.1} (still in the case $k=n$) 
it remains to find a  $\widetilde J$-convex Stein domain $\widetilde W \supset K$ 
contained in $V$ (so that $\wt f$ will be holomorphic in $\wt W$)
and satisfying the other required properties; see Fig.\ \ref{Fig6}.  
Assuming as we may that $\Sigma$ has been standardized along $G(S)$ as described above,
this is an immediate application of Lemma \ref{Eliashberg} --- one takes 
$\wt W =\wt G({\rm Int}L)$ where $L\subset \C^n$ is a standard 
strongly pseudoconvex handlebody around $Q_\lambda \cup D$ as in the cited lemma.
Eliashberg also showed how to extend a $J$-convex defining function 
$\rho$ for $W$ to a $\wt J$-convex defining function $\wt \rho$ for $\wt W$ 
with precisely one additional Morse critical point of index 
$\dim_\R D$ which may be placed at the center of the attached
handle. It follows in particular that $W$ and $\wt W$ are two
sublevel sets of the same $\wt J$-convex exhaustion function and 
hence $W$ is Runge in $\wt W$. This also follows from our 
earlier argument on polynomial convexity of $Q_\lambda\cup D$.  

%
%
%
%

\begin{figure}[ht]
\psset{unit=0.6cm} 
\begin{pspicture}(-4,-4)(9,4)

\pscircle[linestyle=dotted,linewidth=1.2pt](0,0){3}                   
\psarc[linewidth=1.5pt](0,0){3}{66}{294}                                   
\psarc[linewidth=1.5pt](0,0){3}{-32}{32}

\pscircle[linewidth=0.5pt,linestyle=dashed,fillstyle=crosshatch,hatchcolor=yellow](0,0){2.5}    
\pscircle[fillstyle=solid,fillcolor=white,linestyle=none](0,0){0.5}
\rput(0,0){$W$}

\psarc[linestyle=dashed](4.25,0){3}{-145}{145}                        

\psdot[dotsize=3pt](1.8,1.75)                                         
\psdot[dotsize=3pt](1.8,-1.75)

\psarc[linewidth=1.5pt,linecolor=red](4.25,0){2.7}{-118}{118}                       
\psarc[linewidth=1.5pt,linecolor=red](4.25,0){3.3}{-125}{125}                       

\psecurve[linewidth=1.5pt,linecolor=red](6,1)(2.64,1.4)(2.6,2)(3,2.4)(4,2.7)        
\psecurve[linewidth=1.5pt,linecolor=red](6,-1)(2.64,-1.4)(2.6,-2)(3,-2.4)(4,-2.7)   
\psecurve[linewidth=1.5pt,linecolor=red](0,4)(1.2,2.75)(2,2.6)(2.38,2.7)(4,5)       
\psecurve[linewidth=1.5pt,linecolor=red](0,-4)(1.2,-2.75)(2,-2.6)(2.38,-2.7)(4,-5)  

\psline[linewidth=0.3pt]{<-}(3.2,0)(4.6,0)
\psline[linewidth=0.3pt]{->}(5,0.4)(5,2.4)
\psline[linewidth=0.3pt]{<-}(7.25,0)(8.5,0)

\rput(9.6,0){core $M$}

\rput(5,0){$\wt W$}

\end{pspicture}
\caption{The handlebody $\wt W$ around $K=\overline W\cup M$} 
\label{Fig6}

\end{figure}
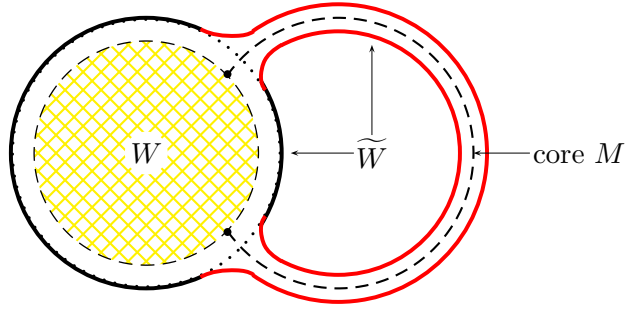

This completes the proof of Lemma \ref{L5.1} for $k=n$. When $1\le k<n$,
we apply the same proof with a totally real $n$-disc $M'$ 
obtained by thickening $M=G(D)$ in the missing $n-k$ real directions. 
To find such $M'$ we choose a $\C$-linear normal framing  
$\beta \colon \nu'' \to TX|_M$ over $G\colon D\hra X$ such that 
$d^\C G\oplus \beta \colon T^\C D\oplus \nu'' = T\C^n|_{D}\to TX|_{G(D)}$ 
is a $\C$-linear isomorphism. Furthermore, we may choose 
$\beta$ to map $\nu''|_S$ into the contact subbundle $\xi|_{G(S)}$.
Let $rD^{n-k}$ denote a closed ball of radius $r>0$ in the real subspace 
$\{0\}^k\oplus \{i0\}^k\oplus \R^{n-k}\oplus \{i0\}^{n-k} \subset \C^n$. 
For a small $r>0$ we can extend $G$ to a smooth $J$-real embedding 
$(1+r) D^k\times r D^{n-k} \hra X$, still denoted $G$, which is real analytic 
near $S^{k-1}\times r D^{n-k}$, it  maps the latter manifold
to a Legendrian submanifold of $\Sigma$, and such that $dG$ 
equals $\beta''$ in the directions tangent to $r D^{n-k}$
at every point of $D^k$. Taking $M'=G(D^k\times r D^{n-k})$ and 
$K' = \overline W \cup M'$ reduces the proof to the case $k=n$.
(The fact that $\di M'$ is not entirely contained in $\di W$
does not cause any complication.)
\end{proof}

%
%
%
%
%
\section{The case $\dim_{\R} X\ne 4$} 
In this section we prove our main results in the case 
when $\dim_\R X\ne 4$. 
Let $P$ be a compact Hausdorff space and let $X$ and $Y$ be smooth manifolds.
A {\em $P$-map} from $X$ to $Y$ is a continuous map 
$f\colon X\times P\to Y$. If $X$ and $Y$ are complex 
manifolds then such $f$ is said to be a {\em holomorphic $P$-map} 
if $f_p =f(\cdotp,p)\colon X\to Y$ is holomorphic for 
every fixed $p\in P$.

\begin{theorem}
\label{T6.1}
Let $(X,J)$ be a smooth almost complex manifold of real dimension 
$2n$ which is exhausted by a Morse function $\rho\colon X\to \R$ 
without critical points of index $>n$. Assume that for some $c\in \R$
the structure $J$ is integrable in $X_c =\{x\in X\colon \rho(x) < c\}$
and $\rho$ is strongly $J$-plurisubharmonic in $X_c$.
Let $Y$ be a complex manifold with a distance function $d_Y$ induced 
by a Riemannian metric. Let $P$ be a compact Hausdorff space 
and $f \colon X\times P\to Y$ be a $P$-map which is 
$J$-holomorphic in $X_c$. 

If $n\ne 2$, or if $n=2$ and $\rho$ has no critical points 
of index $>1$ in $\{x\in X\colon \rho(x)\ge c\}$ 
then for every compact set $K\subset X_c$ 
and for every $\e>0$  there exist a Stein structure $\wt J$ on $X$ 
and a homotopy of $P$-maps $f^t\colon X\times P \to Y$ $(t\in [0,1])$ 
satisfying the following properties:
\begin{itemize}
\item[(a)]  $f^0=f$. 
\item[(b)] The $P$-map $\wt f:=f^1$ is $\wt J$-holomorphic on $X$. 
\item[(c)] There is a homotopy  $J_t$ of almost complex 
structures on $X$ which is fixed in a neighborhood of $K$
such that $J_0=J$ and $J_1=\wt J$; if $J$ is integrable on $X$
then $J_t$ can be chosen integrable for all $t\in[0,1]$.
\item[(d)]  
For every $t\in [0,1]$ the $P$-map $f^t$ is $J$-holomorphic in a 
neighborhood of $K$ and satisfies 
$\sup \{ d_Y(f^t(x,p),f(x,p)) \colon x\in K,\  p\in P\} < \e$.
\end{itemize}
\end{theorem}

The special case of Theorem \ref{T6.1}, 
when applied to the constant map $X\to point$, 
coincides with Theorem 1.3.1 of Eliashberg \cite{E}.

\begin{proof}
We shall present the proof for the nonparametric case since 
the parameters do not present any essential complication. 

Fix a compact set $K\subset X_c$ and choose a regular value 
$c_0\in \R$ of $\rho$ such that $K\subset X_{c_0} \Subset X_c$.
Hence the structure $J=J_0$ is integrable 
in a neighborhood of $\overline X_{c_0}$ and
the map $f\colon X\to Y$ is $J_0$-holomorphic 
in a neighborhood of $\overline X_0$.

Let $p_1,p_2,\ldots$ be the critical points 
of $\rho$ in $\{x \in X \colon \rho(x)>c_0\}$,
ordered so that $\rho(p_j)< \rho(p_{j+1})$ for every $j$.
Choose numbers $c_j$ satisfying 
$c_{-1}=-\infty< c_0 < \rho(p_1)< c_1 < \rho(p_2) < c_2 <\ldots$.
Let $k_j$ denote the Morse index of $p_j$, so $k_j\le n$.
For each $j=0,1,\ldots$ we set $X_j=\{x\in X\colon \rho(x) < c_j\}$,
$\Sigma_j=\di X_{j} =\{x\colon \rho(x)=c_j\}$,
and $A_j = \{x\in X\colon c_{j-1}\le \rho(x) < c_j\}$.
We shall inductively construct a sequence of almost complex structures 
$J_j$ on $X$ and a sequence of maps $f_j\colon X\to Y$ satisfying the following 
for $j=0,1,2,\ldots$:
\begin{itemize}
\item[(i)]  $J_j$ is integrable in a neighborhood of $\overline X_j$
and the manifold $(X_j,J_j)$ is Stein with strongly pseudoconvex boundary,
\item[(ii)] $J_j=J_{j-1}$ in a neighborhood of $\overline X_{j-1}$,
\item[(iii)] the set $\overline X_{j-1}$ is $\cH(X_j,J_j)$-convex, 
\item[(iv)] the map $f_j$ is $J_j$-holomorphic in a  
neighborhood of $\overline X_j$,
\item[(v)]  $\sup_{x\in X_{j-1}} d(f_j(x),f_{j-1}(x)) < \e\, 2^{-j-1}$, and
\item[(vi)] there is a homotopy from $f_{j-1}$ to $f_j$
which is $J_j$-holomorphic and uniformly close to $f_{j-1}$
in a neighborhood of $\overline X_{j-1}$
(satisfying the estimate in (v)).
\end{itemize}

These conditions clearly hold for $j=0$, and in this case  
(ii), (iii) and (v) are vacuous. Assume inductively that the above
hold for $j-1$. By Morse theory $X_j$ is diffeomorphic to a 
handlebody obtained by attaching to $X_{j-1}$ 
an embedded disc $M_j \subset X\bs X_{j-1}$ of dimension $k_j$ 
and smoothly thickening the union $\overline X_{j-1}\cup M_j$ inside $X_j$. 
($M_j$ may be taken as the unstable manifold of the critical point 
$p_j$ for the gradient flow of $\rho$.)  

Applying Lemmas \ref{L3.1} and \ref{L5.1} with $W=X_{j-1}$, $J=J_{j-1}$ and $f=f_{j-1}$ 
we get a Stein structure $\widetilde J$ on a handlebody $\widetilde W$ which is isotopic to 
$X_j$ and satisfies $\overline X_{j-1} \subset \widetilde W \subset X_j$,
and a map $\widetilde f \colon X\to Y$, homotopic to $f_{j-1}$,
which is $\widetilde J$-holomorphic on $\widetilde W$ and
which approximates $f_{j-1}$ uniformly on $X_{j-1}$.
(If $k_j=0$, a new connected component of the sublevel set 
$\{\rho<c\}$ appears at $p_j$ when $c$ passes the value $\rho(p_j)$, 
and it is trivial to find $\wt f$ and $\wt J$ with these properties.)
There is a smooth diffeotopy $h_t\colon X\to X$ $(t\in [0,1])$ 
which is fixed in a neighborhood of $\overline X_{j-1}$ such 
that $h_0$ is the identity map on $X$ and $h=h_1$ satisfies 
$h(X_j)=\widetilde W$. Taking $J_j=h^* (\widetilde J)$ and 
$f_j=\widetilde f\circ h$ completes the inductive step. 
(The homotopy from $f_{j-1}$ to $f_j$ is obtained by composing 
the homotopy from $f_{j-1}$ to $\widetilde f$ by the map $h$.)
The induction may proceed.

By properties (i) and (ii) there is a unique integrable 
complex structure $\widetilde J$ on $X$ which agrees 
with $J_j$ on $X_j$. By the construction $\widetilde J$ is 
homotopic to the initial structure $J=J_0$ since at the 
$j$-th stage of the construction, the structure $J_j$ was chosen 
homotopic to $J_{j-1}$ by a homotopy 
which is fixed on a neighborhood of $\overline X_{j-1}$. 
The complex manifold $(X,\widetilde J)$ is exhausted by 
the increasing sequence of Stein domain $X_j$, and the
Runge property (iii) implies that $(X,\widetilde J)$ is Stein. 
Properties (iv) and (v) insure that the sequence
$f_j \colon X\to Y$ converges uniformly on compacts in $X$ 
to the $\widetilde J$-holomorphic map 
$\widetilde f=\lim_{j\to\infty} f_j\colon X\to Y$
satisfying $\sup_{x\in X_0} d(\widetilde f(x), f_0(x)) < \e$.
By (vi) the homotopies from $f_{j-1}$ to $f_j$ also 
converge, uniformly on compacts in $X$, and give a homotopy 
from the initial map $f_0$ to $\widetilde f$, thus
completing the proof.

If the initial structure $J$ on $X$ is integrable then all
steps can be made within the class of integrable
structures. 
\end{proof}

\smallskip
{\em Proof of Theorem \ref{Main1bis}.} 
This requires only minor modifications of the proof of Theorem \ref{T6.1}.
The main difference is that we do not 
change the given integrable structure $J$ during the construction,
at the  cost of remaining on subsets of $X$ which are only diffeomorphic
to sublevel sets of $\rho$ (and not equal to them as before).
We will in fact obtain a stronger version with approximation,
similar to Theorem \ref{T6.1}. 

We assume the same assumptions and notation
as in the proof of Theorem \ref{T6.1}.
Thus, $W_0=X_0$ is a sublevel set of a Morse exhaustion function
$\rho\colon X\to \R$ which has no critical points of index $>n$
in $X \bs W_0$, $\rho$ is strongly plurisubharmonic 
in a neighborhood of $\overline W_0$, and 
the initial map $f_0=f \colon X\to Y$ is holomorphic in a 
neighborbood of $\overline W_0$.
Let $X_j=\{\rho<c_j\}$ where the constants $c_j$ are chosen as 
in the proof of Theorem \ref{T6.1}, so $\rho$ has a unique critical 
point $p_j$ in $X_j \bs X_{j-1}$. Choose $\e>0$ and let
$d_Y$ denote a distance function on the manifold $Y$.
Assuming that $n\ne 2$ (the case $n=2$ will be treated in \S 7 below) 
we inductively construct an increasing sequence of relatively compact, 
strongly pseudoconvex domains $W_1\subset W_2\subset \cdots \subset X$ 
with smooth boundaries, a sequence of maps $f_j\colon X\to Y$, 
and a sequence of diffeomorphisms $h_j\colon X\to X$ such that the following
hold for all $j=1,2,\ldots$:
\begin{itemize}
\item[(i)] 
$\overline W_{j-1}$ is $\cH(W_j)$-convex, 
\item[(ii)]
$f_j$ is holomorphic in a neighborhood of $\overline W_j$
and is homotopic to $f_{j-1}$ by a homotopy $f_{j,t}\colon X\to Y$
$(t\in [0,1])$ such that each $f_{j,t}$ is holomorphic near 
$\overline W_{j-1}$ and satisfies
$\sup_{x\in W_{j-1}} d_Y\bigl(f_{j,t}(x),f_{j-1}(x)\bigr) < \e 2^{-j}$,
\item[(iii)] 
$h_j(X_j)=W_j$, and 
\item[(iv)] 
$h_j=g_j\circ h_{j-1}$ 
where $g_j\colon X\to X$ is a diffeomorphism of $X$ which is 
diffeotopic to $id_X$ by a diffeotopy which is fixed
in a neighborhood of $\overline W_{j-1}$.
(In particular, $h_j$ agrees with $h_{j-1}$ near $\overline W_{j-1}$.)
\end{itemize}

Granted such sequences, it is easily verified that the limit map 
\[ 
	\wt f=\lim_{j\to \infty} f_j\colon \Omega=\cup_{j=1}^\infty W_j \to Y
\]
and the limit diffeomorphism $h=\lim_{j\to\infty} h_j \colon X\to \Omega$
satisfy the conclusion of Theorem \ref{Main1bis}.

To prove the inductive step we begin by attaching
to $W_{j-1}=h_{j-1}(X_{j-1})$ the disc $M_j := h_{j-1}(D_j)$, 
where $D_j\subset X_j \bs X_{j-1}$ (with $\di D_j\subset \di X_{j-1}$)
is the unstable disc for $\rho$ at the unique critical point 
$p_j \subset X_j \bs X_{j-1}$ of $\rho$ in this region.
By Lemma \ref{L3.1} we can isotope $M_j$ to a totally real, 
real analytic disc in $X$ attached to $\di W_{j-1}$ along a Legendrian
sphere. Applying Lemma  \ref{L5.1} with the integrable structure $J$
we find the next map $f_j\colon X\to Y$ which is holomorphic 
in a thin handlebody $W_j\supset \overline W_{j-1}\cup M_j$. 
The next diffeomorphism $h_j$ with the stated properties is 
then furnished by the Morse theory. 
This concludes the proof of Theorem \ref{Main1bis}.

With a bit more care one can insure that $\di \Omega$ is
smoothly bounded and strongly pseudoconvex, but in general we cannot 
choose such $\Omega$ to be relatively compact, unless $X$ 
admits an exhaustion function $\rho\colon X\to\R$ with at most 
finitely many critical points.
\qed

\smallskip

%
%
%
%
In the remainder of this section we discuss the existence of 
holomorphic maps of maximal rank (immersions resp.\ submersions).
Let $X$ and $Y$ be complex manifolds. A necessary condition for 
a continuous map $f\colon X\to Y$ to be homotopic to a holomorphic 
map of maximal rank is that $f$ is covered by a complex vector bundle
map $\iota\colon TX\to TY$, i.e.,  such that for every 
$x\in X$ the map $\iota_x\colon T_x X\to T_{f(x)} Y$ 
is $\C$-linear and of maximal rank. If $X$ is Stein, this condition 
is known to be also sufficient in the following cases:
\begin{itemize}
\item[(i)]   $\dim X=1$ and $Y=\C$ (Gunning and Narasimhan \cite{GN});
\item[(ii)]  $Y=\C^q$ with $q>\dim X$ 
(Eliashberg and Gromov \cite{EG1,Gbook});
\item[(iii)] $Y=\C^q$ with $q<\dim X$ (Forstneri\v c \cite{FActa});
\item[(iv)]  $n=\dim X\ge \dim Y$ and $Y$ satisfies a 
Runge approximation property for holomorphic submersion  
$\C^n \to Y$ on compact convex sets in $\C^n$ 
(the Property ${\mathrm S}_{\mathrm n}$ in \cite{FFourier1}).
\end{itemize}

By an obvious modification of the proof of Theorem \ref{T6.1},
using the part of Lemma \ref{L5.1} for maps of maximal rank,
one obtains the following result which in particular 
implies Theorem \ref{Main2}. We leave out the details.

\begin{theorem}
\label{T6.3}
Let $(X,J)$ be a smooth almost complex manifold 
of real dimension $2n$, exhausted 
by a Morse function $\rho\colon X\to\R$ without critical points of 
index $>n$. Let $f\colon X\to Y$ be a continuous map to a complex 
manifold $Y$, and let $\iota\colon TX\to TY$ be a complex vector bundle 
map  covering $f$ such that $\iota_x\colon T_x X\to T_{f(x)}Y$ is of 
maximal rank $\min\{\dim X, \dim Y\}$ for every $x\in X$.
If $n\ne 2$, or if $n=2$ and $\rho$ has no critical points 
of index $>1$, there is a homotopy $(J_t,f_t,\iota_t)$ $(t\in [0,1])$ 
where $J_t$ is an almost complex structure on $X$, $f_t\colon X\to Y$
is a continuous map, and $\iota_t\colon TX\to TY$ is a $J_t$-complex 
linear vector bundle map of pointwise maximal rank covering $f_t$, 
such that the following hold:
\begin{itemize}
\item[(i)] $J_0=J$, $f_0=f$, $\iota_0=\iota$, 
\item[(ii)] $(X,J_1)$ is a Stein manifold, 
\item[(iii)] the map $f_1\colon X\to Y$ is $J_1$-holomorphic 
and of maximal rank (an immersion resp.\ a submersion), 
and $df_1=\iota_1$.
\end{itemize}
If in addition there is a constant $c\in \R$ such that $J$ is integrable Stein on 
the set $X_c=\{\rho<c\}$, $f$ is holomorphic on $X_c$ 
and $\iota=df$ on $X_c$ then for every compact set 
$K\subset X_c$ the homotopy $J_t$ may be chosen fixed near $K$, 
the map $f_t$ may be chosen holomorphic near $K$ and uniformly 
close to $f=f_0$ on $K$, and $\iota_t$ may be chosen 
to satisfy $\iota_t=df_t$ near $K$ for each $t$.

The analogous result holds for a family of maps parametrized 
by a compact Hausdorff space (compare with Theorem \ref{T6.1}). 
\end{theorem}

%
%
%
%

\section{The four dimensional case}
The following is a precise version of Theorem \ref{Main1} in 
the case $\dim_\R X=4$. The notion of a $P$-map  
was defined at the beginning of \S 6.

\begin{theorem}
\label{T7.1}
Let $X$ be a smooth oriented 4-manifold, exhausted by a Morse function 
$\rho\colon X\to \R$ without critical points of index $>2$. 
Assume that for some $c\in \R$ there is an integrable
complex structure $J$ on $X_c =\{x\in X\colon \rho(x) < c\}$ 
such that $\rho|_{X_c}$ is strongly $J$-plurisubharmonic.
Let $Y$ be a complex manifold with a distance function 
$d_Y$ induced by a Riemannian metric, let $P$ be a compact 
Hausdorff space, and let $f \colon X\times P\to Y$  be 
a $P$-map which is $J$-holomorphic in $X_c$. 

Given a compact set $K\subset X_c$ and an $\e>0$,
there are a Stein surface $(X',J')$, an orientation 
preserving homeomorphism $h\colon X\to X'$ 
which is biholomorphic in a neighborhood of $K$,
and a holomorphic $P$-map $f'\colon X'\times P \to Y$ 
such that the $P$-map $\wt f\colon X\times P \to Y$,
defined by $\wt f(x,p) = f'(h(x),p)$, is homotopic to $f$
and satisfies
\[
	\sup \left\{ d_Y\bigl( f(x,p), \wt f(x,p) \bigr) 
	\colon x\in K,\ p\in P \right\}  < \e.
\]	
\end{theorem}

Unlike in the case $n>2$, we do not need to assume that the 
almost complex structure $J$ is defined 
on all of $X$ since the obstruction to extending $J$ across 
an attached handle only appears for handles of index $>2$. 
However, if $J$ is already given on all of $X$, one can choose 
$(X',J')$ such that the homotopy class of almost complex structures 
on $X$ determined by $h^*(J')$ equals the class of $J$; Gompf showed 
that this notion makes sense under orientation preserving 
homeomorphisms \cite[p.\ 645]{Go1}.

Before proving Theorem \ref{T7.1} we indicate some consequences.
The following is obtained by combining Theorem \ref{T7.1} 
with Corollary 3.2 and Theorem 3.3 of Gompf \cite[p.\ 648]{Go1}.

\begin{corollary}
\label{Gompf1}
Let $X$ be a smooth, closed, oriented 4-manifold. 
There exists a smooth, finite wedge of circles
$\Gamma\subset X$ such that for every continuous map
$f\colon X\bs \Gamma \to Y$ to a complex manifold $Y$
there is a (possibly exotic) Stein structure on $X\bs \Gamma$ 
and a  holomorphic map $\wt f \colon X\bs \Gamma\to Y$ homotopic to $f$. 
If $X=\C\P^2$, this holds after removing a single point 
(in this case any Stein structure on $\C\P^2\bs\{p\}$ is exotic). 
The analogous result holds for each open oriented 4-manifold 
after removing a suitably chosen smooth 1-complex.
\end{corollary}

The point is that there is a wedge of circles $\Gamma$ in $X$ 
such that $X\bs \Gamma$ admits a handle decomposition without 3-
and 4-handles. The projective plane $\C\P^2$ has a single 4-cell (and no 3-cells)
in its handlebody decomposition, hence removing a point leaves only
cells of index $\le 2$. 

Here is another result obtained by combining Gompf's methods in \cite{Go2} 
with the proof of Theorem \ref{Main1bis} (ii) given below.

%
%
%
%
%
\begin{corollary}
\label{Gompf2}
Let $M$ be a tame, topologically embedded CW 2-complex
in a complex surface $X$ and let $U$ be an open neighborhood of $M$ in $X$.
For every continuous map $f\colon M\to Y$ to a complex manifold $Y$
there exist a topological isotopy $h_t\colon X\to X$,
with $h_0=id_X$ and $h_t(M)\subset U$ for all $t\in [0,1]$,
a Stein thickening $\Omega\subset U$ of the CW complex 
$h_1(M)$, and a holomorphic map $\wt f\colon \Omega\to Y$ 
such that $\wt f\circ h_1 \colon M\to Y$ is 
homotopic to $f$.
\end{corollary}

Gompf showed that the necessary adjustment of the 
initial 2-complex $M$ in $X$ is quite mild from the topological 
point of view, and all essential data of the topological embedding 
$M\hra X$ can be preserved. Stein domains $\Omega$ obtained 
in this way will 
typically have nonsmooth boundaries in $X$ and may be chosen to realize 
uncountably many distinct diffeomorphism types. 
In certain special cases when the 2-cells in $M$ satisfy
certain framing conditions, it is possible to find Stein thickenings 
of a $\cC^0$-small smooth perturbation of $M$ in $X$ which even 
have the {\em diffeomorphism type} of a smooth handlebody with core $M$.
In this direction see also  Costantino \cite{Co}.

\smallskip
{\em Proof of Theorem \ref{T7.1}.}
We shall follow Gompf's construction 
of exotic Stein structures \cite{Go1}, but with a modification 
which will better suit our task of finding a holomorphic map 
in the given homotopy class. Subsequently we will show how the construction 
can be carried out inside a given complex surface as in \cite{Go2}, 
thereby proving Theorem \ref{Main1bis} (for $\dim_\R X=4$) and Corollary 
\ref{Gompf2}. 

The proof of Theorem \ref{T6.1} applies without any changes when attaching
handles of index zero or one, but the difficulty arises 
when attaching 2-handles because Lemma \ref{L3.1} may fail
to give an embedded core 2-disc attached along a Legendrian curve.
As shown by Gompf \cite{Go1,Go2} the obstruction can be removed
by using {\em Casson handles}, but the price to pay is a change 
of the underlying smooth structure on $X$.

We begin by reviewing the necessary background material,
referring to \cite{Go1} for a more complete discussion.
Let $W$ be a relatively compact, smoothly bounded domain 
in $X$ such that $J$ is defined on a neighborhood of $\overline W$
and $\Sigma=\di W$ is strongly $J$-pseudoconvex.
Let $G\colon D=D^2\hookrightarrow X\bs W$ be an
embedded 2-disc attached along the circle $G(S) \subset \Sigma:=\di W$. 
Let $g=G|_S\colon S\hra \Sigma$. The restriction of the contact subbundle 
$\xi=T\Sigma \cap J(T\Sigma)$ to the circle $C:=g(S)$ is a trivial bundle 
(every oriented two-plane bundle over a circle is trivial).
As in \S 3 above we can use the Legendrization theorem 
to make $G$ normal to $\di W$ and its boundary circle 
$C\subset \Sigma$ Legendrian in $\Sigma$.
(The boundary circle remains in $\Sigma$ during this isotopy.)
We denote by $M=G(D)$ the resulting embedded 2-disc 
in $X\bs W$, with $\di M=C$.

Let $\nu_C \subset T\Sigma|_C$ denote the normal bundle of $C$ in $\Sigma$.
It is spanned by the pair of vector fields $(Jw,J\tau)$ 
where $w$ is normal to $\Sigma$ in $X$, with $Jw\in T\Sigma$,
and $\tau$ is tangent to $C$.
The pair $(Jw,J\tau)$ is a {\em canonical framing}, 
or a {\em Thurston-Bennequin framing}, $TB$, of the normal bundle
$\nu_C$. (This notion is only defined for Legendrian knots or links.)

Let $\beta\colon \nu \to TX|_M$ denote a normal framing over $M$
(a trivialization of the normal bundle of $M$ in $X$), chosen such that 
$\beta(\nu|_S) = \nu_C$. We thus have two framings of $\nu_C$,
namely $\beta$ (which extends to the disc $M$)
and the $TB$ framing. Since $\nu_C$ is a trivial 2-plane
bundle over $C$, any two framings differ 
up to homotoy by a map $C\to SO(2)=S^1$, hence by an integer. 
We can thus write $[\beta] = TB  + k$;  the integer
$k=k([\beta])\in\Z$ will be called the {\em framing index} 
of $\beta$. 

In the model case when $M=D^2 \subset \C^2$ is the 
core of a standard handle in $\C^2$ attached to a quadric domain 
$Q_\lambda \subset \C^2$ (\ref{Qlambda}) we easily see that  
\begin{equation}
\label{TB}
		[\beta] = TB  -1.                    
\end{equation}
Indeed, the tangent field $\tau$ to $S=\di D^2$ rotates once 
in the positive (counterclockwise) direction 
as we trace $S$ in the positive direction. 
Since the complex structure operator $J_{st}$ 
on $T_z\C^2$ is an orientation reversing map of $\R^2$ onto 
$i\R^2$, the vector field $J\tau$ (which determines $TB$) rotates once 
in the clockwise direction, hence $\beta$ is obtained from the $TB$ 
framing by one left (negative) twist, so (\ref{TB}) holds. 

When the normal framing $\beta$ of $M$ satisfies (\ref{TB}) then $J$
extends to an integrable complex structure 
in a neighborhood of $\overline W \cup M$ in $X$ such that the 
core disc $M$ is $J$-real (this is precisely as in \cite{E}). 
In this case Lemma \ref{L5.1} in \S 5 applies and yields 
a holomorphic map in a neighborhood of $\overline W\cup M$
which approximates the previous map uniformly on $W$. 
If this ideal situation occurs for all 2-handles in $X\bs W$ 
then the construction of a Stein structure on $X$, and of a
holomorphic map $X\to Y$ in a given homotopy class,
can be completed exactly as in \S 6.

%
%

Suppose now that $k=[\beta]-TB  \ne -1$ for some 2-handle. 
A basic fact from the theory of Legendrian knots \cite{Ar,E2} 
is that for any Legendrian knot $K$ there is a $\cC^0$-small isotopy 
preserving the knot type, but changing its Legendrian knot type, which 
adds a desired number of left (negative) twists to the $TB$ framing. 
(One adds small spirals to $K$.) Since the homotopy class of the $\beta$ 
framing is preserved under an isotopy of $C$ in $\Sigma$, 
we see that $k=[\beta]-TB$ can be increased by any number of units.
If $k<-1$, it is therefore possible to add spirals to the boundary circle 
and obtain an isotopic embedding $(D,S) \hookrightarrow (X\bs W,\Sigma)$  
satisfying (\ref{TB}), thereby reducing the problem to the previous case. 

%
%
%
%
The problem is more difficult when $k\ge 0$ since it is in 
general impossible to add right twists to the $TB$ framing 
(equivalently, to decrease the framing index $k$). 
This is only possible in a contact structure 
which is {\em overtwisted}, in the sense that it contains
a topologically unknotted Legendrian knot $K$ with the Thurston-Bennequin
index $tb(K)=0$; adding such knot to a Legendrian knot adds a 
positive twist to the $TB$ framing, making it possible to 
decrease $k=[\beta]-TB$ and hence reach $k=-1$. 
However, Eliashberg proved in \cite{E2} that contact structures 
arising as boundaries of strongly pseudoconvex Stein manifolds 
are never overtwisted (they are {\em tight}). A $2$-handle  
for which we cannot find an isotopy of the boundary circle 
to a Legendrian knot so that (\ref{TB}) holds will be called 
in the sequel {\em a wrongly attached handle}.

In \cite{Go1} Gompf showed how one can circumvent the 
problem by replacing each wrongly attached $2$-handle 
by a {\em Casson handle} which is homeomorphic,
but not diffeomorphic, to the standard 2-handle $D^2\times D^2$. 
In such case the new domain can be chosen to admit a 
Stein structure and is homeomorphic, but not diffeomorphic, 
to the original manifold. In the following two paragraphs we 
summarize Gompf's construction for future reference. 

Let $Z$ be a smooth $4$-manifold obtained by attaching 
a $2$-handle $h$, with the core disc $M$, to 
a compact Stein surface $\overline W$ along its 
strongly pseudoconvex boundary $\Sigma=\di W$. We first isotope 
$C=\di M\subset \Sigma$ inside $\Sigma$ to a Legendrian knot. Since the $TB$ invariant 
can be increased by a non-Legendrian isotopy by an arbitrary integer,
while the homotopy class $[\beta]$ does not change by the isotopy,  
we can assume that the framing coefficient $[\beta]-TB$ is odd,
$[\beta]-TB =-1 +2k$ for some $k\in \Z$. If $k<0$, 
we can isotope the boundary of $h$ (by adding left twists
to the $TB$ framing) to get (\ref{TB}) and we are done. 
If not, we remove the 2-handle 
$h$ and reattach it to $W$ along $C$ using the framing 
$TB-1$, meaning that we add $2k$ negative twists to the 
framing of $h$. Lets us call this new handle $h'$,
and let $Z'$ be the new manifold obtained by attaching 
$h'$ to $\overline W$ in this way.
As mentioned above, we can extend the Stein structure 
from $\overline W$ across the 2-handle $h'$ to a Stein structure on $Z'$. 
The problem is that the manifold $Z'$ is not diffeomorphic to $Z$.
In fact, if the 2-handle $h$ in the manifold $Z$ gives a homology class in 
$H_2(Z,\Z)$ with self-intersection $m$, the self-intersection 
of the handle $h'$ in the homology of $Z'$ equals $m-2k$. In order 
to compensate this error we now make $k$ {\em positive self-plumbings} 
to the handle $h' \subset Z'$ to get a yet new manifold $Z_1$. 
(A positive self-plumbing is done by choosing two disjoint 
closed discs $\disc_1,\disc_2$ in the core disc $M$ of the handle $h'$, 
trivializing the normal bundle of $M$  over these discs
to get subsets $\disc_1\times \disc'_1$, $\disc_2\times \disc'_2$ 
of $Z'$ diffeomorphic to the standard handle $D^2\times D^2$,  
and identifying them by the map  $\disc_1\times \disc'_1 \to \disc_2\times \disc'_2$ 
that interchanges the factors: $(z,w)\mapsto (w,z)$ in the coordinates on $D^2\times D^2$. 
(For a negative plumbing we would use $(z,w)\to (\overline w,z)$.) 
The manifold $Z_1$ is the quotient space of $Z'$ under the above 
identification. The image of $M$ is an immersed 2-disc $M_1$ in $Z_1$ with exactly 
one positive transverse double point $p$ (the image of the centers $p_1$ resp.\ $p_2$
of the discs $\disc_1$ resp.\ $\disc_2$). 
The self-plumbing can be done so that in local coordinates near $p$,
$M_1$ equals $\R^2\cup i\R^2\subset \C^2$,
with $p$ corresponding to the origin, and $Z_1$ is a tubular neighborhood of this set.
Since $\R^2\cup i\R^2$ admits a basis of tubular Stein neighborhoods in $\C^2$
(see \cite[Theorem 1.3.5.]{E} or \cite[Theorem 2.2]{FSteindomains}),
we get a Stein structure on $Z_1$ which agrees 
with the original structure inherited from $Z'$ away
from the plumbed double point $p$. 
Although $Z_1$ is not even homotopically equivalent to $Z$, 
we do have $H_2(Z,\Z)$ isomorphic to $H_2(Z_1,\Z)$ 
by an isomorphism preserving the self-intersection form. 
This follows from the fact that a positive self-plumbing 
of a handle $h'$ introduces a positive transverse double point 
to its core, thus raising the self-intersection number of $h'$ by $2$. 

In the next step we show that $Z$ can be reconstructed 
back from the manifold $Z_1$ by attaching additional 2-cells
to $Z_1$. The group $H_1(Z_1,\Z)$ differs from $H_1(Z,\Z)$ by $\Z^{k}$, 
with new homology classes represented by one loop 
in $Z_1$ for each performed plumbing, the only requirement 
being that the loop passes once through 
the plumbed double point. Attaching a 2-handle to $Z_1$ along 
each such loop cancels the extra homology and moreover, 
if the framing is correct, reconstructs the original manifold $Z$. The details can 
be found in \cite{Go1}, but will also be obvious from the construction below. 
The problem now is that the framing of the 2-handles which we need 
to attach to $Z_1$ (in order to reconstruct $Z$) may not be correct, 
in the sense that the attaching circles of the core discs cannot be isotoped 
to Legendrian knots satisfying (\ref{TB}). To correct this, 
one repeats the above steps, adding kinks to each of these 
wrongly attached handles, thus beginning the {\em Casson tower}
procedure. In this way one gets an increasing sequence of 
Stein manifolds $X_1\subset X_2\subset X_3\subset \cdots$ with 
strongly pseudoconvex boundaries, each of them Runge in the next one.
The limit manifold $\cup_j X_j$ is then also Stein, and by Freedman's result 
on Casson handles \cite{Freedman} it is homeomorphic to the 
original manifold~$X$.

%
%
%
%

Here is a somewhat different explanation of the above procedure
which is better suited to our purpose; its main advantage
is that we remain inside the same manifold $X$ during 
the entire construction.

Let us take an immersed totally real sphere $S^2\to \C^2$ with 
a positive double point at $0\in\C^2$ and no other double points; 
an explicit Lagrangian example is due to Weinstein \cite{Wein1}: 
\begin{equation}
\label{Wkink}
	F(x,y,z) = \bigl(x(1+2iz),y(1+2iz)\bigr)\in \C^2
\end{equation}
where $(x,y,z)\in\R^3$, $x^2+y^2+z^2=1$.
Let us think of $S^2$ as the union of two closed discs 
$D_0\cup D_\infty$, glued along their boundary 
circles $S^1=D_0\cap D_\infty$ and chosen such that 
$F(D_\infty) \subset \C^2$ is embedded while $F(D_0)$ 
contains the positive double point at the origin.
The oriented normal bundle $\nu$ of $F$
is isomorphic to $\overline {TS^2}$, the tangent bundle
of $S^2$ with the reversed orientation. Indeed, 
$T\C^2|_{F(S^2)} = F_*(TS^2)\oplus \nu$,
and the complex structure $J_{st}$ gives an orientation 
reversing isomorphism of the first onto the second 
summand. By a small modification of  $F$ we may assume that 
the double point at $0$ (the {\em center of the kink}) 
is special (\ref{specialdouble}), meaning that in a suitable local 
holomorphic coordinate system the image equals $\R^2\cup i\R^2$
near the origin. We shall take $K=F(D_0) \subset \C^2$ 
as our {\em standard kink} which will be used to correct the 
framing coefficient of wrongly attached handles.

\begin{figure}[ht]
\psset{unit=0.7cm, linewidth=0.7pt} 

\begin{pspicture}(0,-3)(16,3)
\definecolor{myblue}{rgb}{0.66,0.78,1.00}

%
%
\psellipse[fillstyle=solid,fillcolor=myblue](4,0)(3,2)
\pscircle[fillstyle=solid,fillcolor=white,linestyle=none](4,0){0.4}
\rput(4,0){$W$}

%
%
\psarc(8,0){2}{45}{143}
\psarc(8,0){2}{-143}{-45}
\psarc[linestyle=dotted,linewidth=1.2pt](8,0){2}{-45}{45}
\psline(9.42,1.42)(10.8,0)   
\psline(9.42,-1.42)(10.8,0)

%
%
\pscustom[fillstyle=crosshatch,hatchcolor=yellow]
{
\psline(10.8,0)(11.6,-0.8)
\psarc(12.4,0){1.13}{-135}{135}
\psline(11.6,0.8)(10.8,0)
}

\psdots(10.8,0)(6.4,1.2)(6.4,-1.2)

\rput(8,1){$M$}
\psline[linewidth=0.2pt]{->}(8,1.3)(8,1.97)
\rput(12.4,0){$\Delta$}

\psline[linewidth=0.2pt]{<-}(10.4,-0.45)(10.7,-1.4)
\psline[linewidth=0.2pt]{->}(10.8,-1.4)(11.4,-0.6)
\rput(10.8,-1.8){kink}
\rput(10.8,0.5){$p$}

\end{pspicture}
\caption{A kinky disc $M$ with a trivializing 2-cell $\Delta$}
\label{Fig7}
\end{figure}
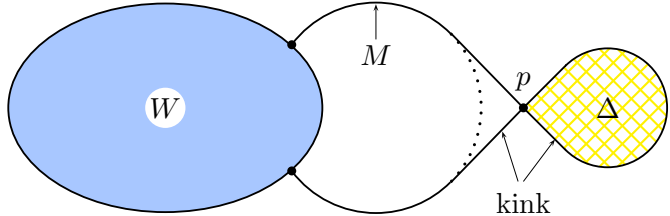

Since $TS$ has the Euler number $\chi(TS)=2$, the normal bundle 
$\nu$ of the Weinstein sphere in $\C^2$ has $\chi(\nu)=-2$; 
hence a copy $K$, glued into a 2-disc $M$
attached to $\di W$ along a Legendian knot as in Fig.\ \ref{Fig7}, 
will reduce the framing coefficient of $M$ by two units. 
This is seen explicitly as follows.
Thinking of $S^2$ as the Riemann sphere $\C\cup\{\infty\}$,
with $D_0=\{z\in\C\colon |z|\le 1\}$ and 
$D_\infty=\{w\in \C\cup\{\infty\} \colon |w|\ge 1\}$,
the real and the imaginary part of the complex vector field 
$\frac{\di}{\di w}$ provide a reference framing for $TS|_{D_\infty}$.
From $\frac{\di}{\di w} = -z^2 \frac{\di}{\di z}$ we see that  
$\frac{\di}{\di w}$ makes two left (negative) twists 
when compared to the framing $\frac{\di}{\di z}$ for $TS^2|_{D_0}$
as we trace the circle $S=\di D_0=\{|z|=1\}$ in the positive direction.
Conversely, $\frac{\di}{\di z}$ makes two right (positive) twists
in comparison to $\frac{\di}{\di w}$. 
Considering these framings on the immersed sphere $F(S^2)\subset \C^2$ 
and applying $J_{st}$ we obtain two framings of the normal bundle 
$\nu$ over the respective discs. Due to the reversal 
of the orientation under $J_{st}$ we see that the framing 
for $\nu|_{D_0}$ makes two left (negative) twists 
when compared to the framing for $\nu|_{D_\infty}$,
which explains $\chi(\nu)=-2$. 

With $F$ as in (\ref{Wkink}), let 
$\Delta= \{F(0,y,z)\colon y\ge 0,\ y^2+z^2\le 1\} \subset \C^2$.
This 2-disc is embedded in $\C^2$, except along the side
$\{y=0\}$ which gets pinched to $0\in \C^2$.
Note that $\di \Delta \subset F(S^2)$, and the union 
$F(S^2) \cup \Delta$ has a tubular neighborhood
diffeomorphic to $S^2\times \R^2$.

In order to make a self-intersection at a point $p$ in the core disc 
$M$ of a handle $h$  in our 4-manifold $X$, we replace a small disc 
in $M$ around $p$ by a copy of the standard kink $K$. 
(See Fig.\ \ref{Fig7}; we removed the small dotted disc and 
smoothly attached along its boundary the kinky disc shown on the 
right.)  We have seen that this surgery reduces the relative Euler number 
over the immersed disc $M$ by $2$ for each kink. 
Adding $k$ kinks on $M$ inside $X$ and then taking a 
tubular neighborhood has the same effect as first 
removing the handle $h$ from $X$, reattaching it with 
a framing of the boundary reduced by $2k$, and then performing 
$k$ self-plumbings on $h$ (as was done by Gompf \cite{Go1}
and described above). In this way we see that the manifold $Z_1$, constructed 
above when discussing Gompf's proof, can be seen as a submanifold of 
the original manifold $X$, changed only by a surgery 
in a small coordinate neighborhood of each of the
kinked points on the core disc of the handle $h$. 
We can also explicitly see the trivializing 2-cell $\Delta$ that needs 
to be added to each of the kinks in order to reconstruct the desired
manifold.

%
%
%
%

In the next stage of the construction every such disc $\Delta$ 
will also have to receive a kink in order to correct its framing coefficient.
This begins the {\em Casson tower} procedure which will converge 
to a Casson handle in place of the original removed disc in $M$.

%
%
%
%

We are now ready to complete the proof of Theorem \ref{T7.1}. 
Assume that our 4-manifold $X$ is constructed by successively 
attaching handles $h_1,h_2,h_3,\ldots$ of index $\le 2$,
beginning with the compact domain $\overline W \subset X$ with smooth boundary 
$\Sigma=\di W$. By assumption we also have an integrable complex 
structure $J$ in a neighborhood of $\overline W$ such that 
$W$ is Stein and its boundary $\Sigma$ is strongly pseudoconvex. 
Let $M_1,M_2,\ldots$ be the cores of the handles $h_1,h_2,\ldots$, 
chosen such that their union is a smoothly embedded $CW$ complex 
inside $X$. Since we have not assumed that our handlebody 
is finite, we can not ask for the ordering of the handles 
with regards to their indices. However, due to local compactness 
we can, and will, ask that when a handle $h_j$ with core $M_j$ is being 
attached, all handles whose core discs intersect the boundary $\di M_j$ 
have already been attached in earlier steps. We can also assume 
that $\di M_j$ consists only of the core discs of handles 
of lower indices. We can now proceed with the induction as in the proof 
of Theorem \ref{T6.1}, but with the following modifications:
\begin{itemize}
\item[(1)] When a 2-handle is attached with a wrong framing, 
we insert the right number of kinks to its core disc
(inside $X$) in order to change the framing coefficient 
to $-1$, thereby insuring that we can extend $J$ to a Stein structure 
in a tubular strongly pseudoconvex neighborhood of the immersed disc. 
(The disc is totally real in this structure, with a special double point 
(\ref{specialdouble}) at each kink.) 
\item[(2)] Each time before proceeding to the next handle $h_{j+1}$, 
we perform one more step on each of the kinked discs 
appearing in the sequence before. More precisely, we add a new 
kinked disc which cancels the superfluous loop at the self-intersection
point introduced in the previous step. (Of course this new kinked
disc introduces a new superfluous loop which will have to be cancelled 
in the subsequent step.)  
\end{itemize}

The first condition is essential since we need to build a 
manifold that is Stein. The second condition insures 
that each handle is properly worked upon, thereby producing a
Casson tower at every place where a kink was made in the
initial 2-disc. At every step we also apply Lemma \ref{L5.1}
to approximate the given map, which has already been made
holomorphic in a tubular strongly pseudoconvex neighborhood 
of our partial (finite) subcomplex, by a map holomorphic 
in a tubular neighborhood of the previous domain with all 
core discs that have been added at the given step.

The proof can now be concluded as in Theorem \ref{T6.1}. 
We construct an increasing sequence of Stein domains  
$X_1\subset X_2\subset \cdots$ inside the original smooth 4-manifold $X$, 
each of them Runge in the next one, together with a sequence of maps 
$f_j\colon X\to Y$ $(j=1,2,\ldots)$ such that $f_j$ is holomorphic 
on $X_j$, it approximates $f_{j-1}$ uniformly on $X_{j-1}$,
and is homotopic to $f_{j-1}$ by a homotopy which is 
holomorphic and uniformly close to $f_{j-1}$ on $X_{j-1}$. 
The Runge property insures that the limit manifold 
$X' =\cup_j X_j$  is Stein with respect to the limit 
complex structure and, by the construction, it is homeomorphic to $X$. 
(It is diffeomorphic to $X$ if no Casson handles were used in the construction.)
A small ambient topological deformation moves the initial 
CW complex (made of cores of the attached handles) 
into $X'$; see \cite{Go2} for more details.  
By construction the limit map 
$f' =\lim_{j\to\infty} f_j \colon X'\to Y$ 
is holomorphic, and the map $f'\circ h\colon X\to Y$ 
is homotopic to $f$.   

The same proof applies to any smoothly embedded 2-complex 
$M$ inside $X$: After a small ambient topological deformation
we find a new embedding $M'\hra X$ with a Stein thickening 
$X'\subset X$ such that a given continuous map $M\to Y$ 
admits a holomorphic representative $X' \to Y$.
\qed

\begin{remark} 
Whenever a handle is wrongly attached, 
the above process is never finite. The reason is 
that in the standard kink $K$, the disc $\Delta$ needed 
to be added to reconstruct the original manifold requires 
exactly one positive kink in order to be able to extend the 
Stein structure to its neighborhood. 
\end{remark}

\begin{remark} 
The complex structure in the above proof was built by 
stepwise extension across handles. If one begins with an almost 
complex structure $J$ defined on all of $X$, 
the Stein structure constructed by this process is in general 
not homotopic to $J$ (compare with Proposition \ref{P3.4}). 
Below we shall give a different construction which will 
produce a Stein structure homotopic to the given initial 
complex structure on $X$.
With a bit more work as in \cite{Go1} one can also show that 
every relative cohomology class in 
$H^2(X,W;\Z)$ which reduces mod $2$ to the second Stiefel-Whitney 
class $w_2(TX)\in H^2(X,W;\Z_2)$ can be obtained 
as a reduction of a cohomology class of $c_1(X,J)$ for an 
appropriate $J$ in Theorem \ref{T7.1}.   
\end{remark}

{\em Proof of Theorem \ref{Main1bis} for $n=2$.} 
We assume that $(X,J)$ is a complex surface
with an integrable (not necessarily Stein) structure $J$,
and with a correct handlebody structure.  
By a modification of the proof given above we shall construct an 
increasing sequence of domains $X_1\subset X_2\subset \ldots$
in $X$ and $J$-holomorphic maps $f_j\colon X_j\to Y$ 
such that each $X_j$ is a $J$-Stein domain in $X$ 
and the other properties are as before. 
The proof which we shall give is similar to the construction 
of Gompf \cite{G2}, the difference being that we 
do everything by using the special kink $K$ introduced above. 

We need to recall a few basic facts about the complex points of 
real surfaces in complex surfaces. Let $M$ be a closed
real surface smoothly immersed in a complex surface $X$. 
A generic such $M$ has finitely many complex points $p\in M$ 
(i.e., points where the tangent space $T_p M$ is a complex line in $T_p X$)
and finitely many transverse double points;  in addition,
no complex point is a double point of $M$.
Locally near a complex point $p$ the surface $M$ is given in suitable
local holomorphic coordinates $(z,w)$, with $z(p)=w(p)=0$,
by an equation $w=|z|^2+\lambda(z^2+\bar z^2)+o(|z|^2)$ for a unique
$\lambda\ge 0$ (Bishop \cite{Bishop}). The point $p$ is {\em elliptic}
if $\lambda<1/2$ and {\em hyperbolic} if $\lambda>1/2$;
the degenerate case $\lambda=1/2$ does not arise for a generic $M$.
At an elliptic point $M$ has a nontrivial local envelope of holomorphy
consisting of a family of small analytic discs (the so called 
Bishop discs, \cite{Bishop}). At a hyperbolic point
$M$ is locally holomorphically convex \cite{FSt} and it admits a basis
of tubular Stein neighborhoods \cite{Sl}.  

If $M$ is oriented, one further divides its complex points  
into positive and negative ones, depending on whether the 
standard complex line orientation of $T_p M$ agrees or disagrees
with the chosen orientation of $M$. 
Let $e_\pm(M)$ resp.\ $h_\pm( M)$ indicate the numbers of 
elliptic resp.\ hyperbolic points in each orientation class.
The {\em Lai indices} $I_\pm( M)=e_\pm( M)-h_\pm( M)$ 
are invariant under a regular homotopy of the immersion 
$M\to X$ and satisfy
\begin{equation}
\label{Lai}
    2\, I_\pm( M)= \chi( M) \pm \langle c_1(X),[M]\rangle + [M]^2 -  2\, d( M).
\end{equation}
Here $\chi(M)$ is the Euler number of $M$, $[M]\in H_2(X;\Z)$ 
is the homology class of $M$, $c_1(X)=c_1(TX)$ 
is the first Chern class of the tangent bundle $TX$,
and $d(M)$ is the algebraic number of self-intersection 
points  of $M$ counted with oriented intersection indices
(see \cite{Bishop}, \cite{CS}, \cite{EH}, \cite{FSteindomains},
\cite{Lai}, \cite{N}, \cite{Web}). 

A basic fact proved by Eliashberg and Harlamov \cite{EH}
(see also \cite{Fcomplexpoints} and \cite{N}) is that
one can cancel a pair of an elliptic and a hyperbolic point
in the same orientation class by a $\cC^0$-small isotopy supported 
in a neighborhood of a suitably chosen embedded arc in $M$ 
connecting these two points. If $I_\pm(M)\le 0$, we can cancel 
all elliptic points and deform $M$ to an immersed surface 
with only special double points (\ref{specialdouble})
which is either totally real 
(if $I_\pm(M)=0$) or has only {\em special hyperbolic points}, 
given in local holomorphic coordinates by $w=z^2+ \bar z^2$
\cite[Theorem 1.2, p.\ 82]{FSteindomains}. 
By \cite[Theorem 2.2]{FSteindomains} this new $M$ admits a  
basis of smoothly bounded Stein neighborhoods diffeomorphic
to $M\times \R^2$, given as sublevel sets of a smooth 
plurisubharmonic function $\tau\ge 0$ which vanishes precisely 
on $M$ and has no critical points in a deleted neighborhood of $M$. 
(Locally at a special hyperbolic point $w=z^2+\bar z^2$ we can take 
$\tau(z,w)=|w-z^2 - \bar z|^2$. 
It was shown in \cite{Sl} how to find strongly pseudoconvex 
tubular neighborhoods of any surface with only hyperbolic complex points.)

By adding a positive kink to $M$ we do not change its homology class 
$[M]$ and its Euler number $\chi(M)$, but we increase by one the number
of double points $d(M)$. Hence we see from (\ref{Lai}) that a positive kink 
reduces each of the numbers $I_\pm( M)$ by $1$. 
This can also be seen directly: 
If in local coordinates $(z,w)=(x+iy,u+iv)$ 
the double point is given by $R^2\cup i\R^2=\{(x+iu)(y+iv)=0\}$
(with the orientation reversed on the second summand),
the perturbation $(x+iu)(y+iv)=2\epsilon^2$ for $\e>0$
is a smooth surface with two hyperbolic complex points 
$p_1=\e(1+i,-1+i)$, $p_2=-p_1$. 
Adding a sufficient number of positive kinks  
and performing a $\cC^0$-small regular homotopy as above
we find an immersed surface $M$ with a basis of strongly 
pseudoconvex tubular Stein neighborhoods. For this argument 
we do not need the formula (\ref{Lai}), 
it suffices to know that a positive kink reduces each of 
the numbers $I_\pm$ by one so that eventually they become nonpositive.

The same argument holds for a compact surface $M\subset X\bs W$ 
with boundary attached to $\Sigma=\di W$ along a 
Legendrian curve --- after kinking it enough times we can find an 
isotopy,  fixed near $\di M$, to a new surface with only special 
hyperbolic points and special double points. (If $M$ is a disc, 
it can be seen that $I_+(M)+I_-(M) = [\beta]-TB+1$, 
but this will not be used.) 

We also note that a continuous function on $M$ can be 
uniformly approximated by functions holomorphic in a
neighborhood of any hyperbolic complex point $p\in M$. 
If $M$ only has hyperbolic complex points then we can 
approximate each continuous function on $M$ uniformly by
holomorphic functions in a Stein neighborhood of $M$.
It is now easily seen that Theorem \ref{T4.1}
still holds (with uniform approximation on $M$) in the presence 
of hyperbolic points. The existence of Stein
neighborhoods  of the graph of $f$ in $X\times Y$, 
needed in the proof of Theorem \ref{T4.1}, is easily insured
by chosing $f$ to be constant in a small neighborhood of 
each hyperbolic point.

We can now complete the proof of Theorem \ref{Main1bis} when $\dim_\R X=4$.
Assume that $X$ is obtained from a strongly $J$-pseudoconvex domain 
$\overline W \subset X$ by successively attaching handles 
$h_1,h_2,\ldots$ with core discs $M_1,M_2,\ldots$,
where the ordering of these handles satisfies the same condition as 
in the proof of Theorem \ref{T7.1}. (We may begin with $W=\emptyset$.)
We shall use the same notation as in the proof of Theorem \ref{Main1bis} 
for $\dim_\R X\ne 4$ (see \S 6), beginning with $W_0=W$ and $f_0=f$. 
The complex structure on $X$ will not change during the proof.

In the inductive step we have a smoothly bounded, 
strongly pseudoconvex domain $W_j\subset X$ 
and a map $f_j\colon X\to Y$ which is $J$-holomorphic in a neighborhod
of $\overline W_j$. (The set $W_j$ is a tubular neighborhood of the union
of $\overline W$ with the cores of handles attached in the earlier steps;
since these cores may have received kinks, 
$W_j$ does not have the correct homeomorphic type, but this will
be corrected in the limit by the Casson handles resulting from the construction.) 
We now attach to $\overline W_j$ the next handle in the sequence. 
As before, attaching a $1$-handle does not pose a problem. 
For a 2-handle $h$ we first make sure that the boundary 
of its core $M$ is a Legendrian curve in $\di W_j$,
and then we add enough positive kinks to algebraically cancel off 
all elliptic points on the core disc $M$ as explained above. 
We denote this new immersed disc by $M'$.  After a $\cC^0$-small regular 
homotopy fixing the boundary $\di M'$ we can assume that $\overline W_j \cup M'$ 
has a basis of tubular, strongly pseudoconvex Stein neighborhoods 
in $X$. We also add to $\overline W_j$ a new trivializing kinky disc 
at each of the kinks from the earlier stages of 
the construction, making sure that the conditions (1) and (2) 
in the proof of Theorem \ref{T7.1} are satisfied.
These additional kinky discs $\Delta_1,\ldots, \Delta_k$ can be chosen 
such that $L_j:=\overline W_j \cup M' \cup (\cup_{l=1}^k \Delta_l)$ 
admits a basis of tubular, strongly pseudoconvex, Stein neighborhoods 
in $X$. Using Lemma \ref{L5.1} we approximate $f_j$ 
uniformly on $L_j$ by a map $f_{j+1} \colon X\to Y$ 
which is holomorphic in a neighborhood of $L_j$.
Choosing a small strongly pseudoconvex tubular neighborhood 
$W_{j+1} \supset L_j$ contained in the domain of holomorphicity
of $f_{j+1}$ completes the induction step.
In the limit we obtain a holomorphic map $f'\colon \Omega\to Y$ 
on the Stein domain $\Omega=\cup_j W_j \subset X$ with the stated properties.
\qed

\smallskip
Corollary \ref{Gompf2} is proved by applying the construction, 
described in the proof of Theorem \ref{Main1bis}, 
to an embedded CW 2-complex in $X$.

\begin{remark}
By a small addition to the above argument it is possible
to change every real surface $M$, attached from the outside
to a strongly pseudoconvex domain $W\subset X$ along a Legendrian link, 
to a totally real immersed surface all of whose double points are  
positive kinks. Indeed, by an isotopy of $\di M\hra \di W$ to another 
Legendrian link one can insure that $I_+(M)-I_-(M)=0$. 
(This difference equals the {\em rotation number} of the canonical 
contact extension of the Legendrian link $\di M \subset \di W$ 
which can be made equal to any given number by isotopies; 
see \cite{E} and \cite{Go1}.) 
If $I_+(M)=I_-(M)>0$, we add this many positive kinks 
and reduce $I_\pm(M)$ to $0$, so our disc becomes totally
real after a small isotopy. If $I_+(M)=I_-(M)<0$, we can increase 
them to zero by adding left spirals to $\di M$ \cite{Go1},
again making $M$ totally real. In this way one can avoid having to deal 
with hyperbolic complex points at the cost of using Legendrian knot theory.
\end{remark}

%
%
%
%
%
\smallskip 
\textit{Acknowledgements.} 
We thank Yakov Eliashberg whose help was essential 
in proving Lemma \ref{L3.1} for handles of maximal dimension,
Robert Gompf for very helpful communication regarding his works 
\cite{Go1,Go2} on Stein surfaces, 
Petar Pave\v si\'c and Sa\v so Strle for their help with 
Proposition \ref{P3.4}, and Irena Majcen for introducing us to
drawing pictures with {\em pstricks}.

\bibliographystyle{amsplain}

\end{document}